\documentclass[11pt]{amsart}
\usepackage{amssymb}

\date{\today}

\newtheorem{teo}{Theorem}[section]
\newtheorem{prop}[teo]{Proposition}
\newtheorem{cor}[teo]{Corollary}

\newtheorem{defin}[teo]{Definition}
\newtheorem{rem}[teo]{Remark}
\newtheorem{ej}[teo]{Example}

\newtheorem{apart}[teo]{ }

\newtheorem{notac}[teo]{Notation}

\newcommand{\ot}{\otimes}
\newcommand{\co}{\circ}

\title[Cohomology of algebras over weak  Hopf algebras]
{Cohomology of algebras over weak  Hopf algebras}

\author[J.N. Alonso]{J.N. Alonso \'Alvarez}
\author[J.M. Fern\'andez]{J.M. Fern\'andez Vilaboa}
\author[R. Gonz\'{a}lez]{R. Gonz\'{a}lez Rodr\'{\i}guez}

\address[J.N. Alonso]{Departamento de Matem\'{a}ticas, Universidad
de Vigo, Campus Universitario Lagoas-Marcosende, E-36280 Vigo,
Spain.} \email[J.N. Alonso]{jnalonso@uvigo.es}

\address[J.M. Fern\'andez]{Departamento de \'Alxebra, Universidad de
Santiago de Compostela.  E-15771 Santiago de Compostela, Spain.}
\email[J.M. Fern\'andez]{alvila@usc.es}

\address[R. Gonz\'{a}lez]{Departamento de Matem\'{a}tica Aplicada II,
Universidad de Vigo, Campus Universitario Lagoas-Marcosende, E-36280
Vigo, Spain.} \email[R. Gonz\'{a}lez]{rgon@dma.uvigo.es}



\begin{document}

\begin{abstract} In this paper we present the Sweedler cohomology
for a cocommutative weak Hopf algebra $H$. We show that the second
cohomology group classifies completely  weak crossed products,
having a common preunit, of $H$ with a commutative left $H$-module
algebra $A$.
\end{abstract}

\maketitle

\section*{Introduction}

In \cite{Moss} Sweedler introduced the cohomology of a cocommutative
Hopf algebra $H$ with coefficients in a commutative $H$-module
algebra $A$. We will denote it as Sweedler cohomology
$H_{\varphi_{A}}(H^{\bullet},A)$ where $\varphi_{A}$ is a fixed
action of $H$ over $A$. If $H$ is the group algebra $kG$ of a group
$G$ and $A$ is an admissible $kG$-module, the  Sweedler cohomology
$H_{\varphi_{A}}^{\bullet}(kG,A)$ is canonically isomorphic to the
group cohomology of $G$ in the multiplicative group of invertible
elements of $A$. If $H$ is the enveloping algebra $UL$ of a Lie
algebra $L$, for $i>1$, the Sweedler cohomology
$H_{\varphi_{A}}^{i}(UL,A)$ is canonically isomorphic to the Lie
cohomology of $L$ in the underlying vector space of $A$. Also, in
\cite{Moss} we can find an interesting interpretation of
$H_{\varphi_{A}}^2(H,A)$ in terms of extensions: This cohomology
group classifies the group of equivalence classes of cleft
extensions, i.e., classes of equivalent crossed products determined
by a 2-cocycle. This result was extended by Doi \cite{doi} proving
that, in the non commutative case, there exists a bijection between
the isomorphism classes of $H$-cleft extensions $B$ of $A$ and
equivalence classes of crossed systems for $H$ over $A$. If $H$ is
cocommutative the equivalence is described by
$H_{\varphi_{\mathcal{Z}(A)}}^{2}(H,\mathcal{Z}(A))$ where
$\mathcal{Z}(A)$ is the center of $A$. Subsequently,  Schauenburg in
\cite{SCH} extended the cohomological results about extensions
including a theory of abstract kernels and their obstructions. The
dual Sweedler theory was investigated by Doi and Takeuchi in
\cite{doi-take} giving a general formulation of a cohomology for
comodule coalgebras for a commutative Hopf algebra as for example
the coordinate ring of an affine algebraic group.

With the recent arise of weak Hopf algebras, introduced by B\"ohm,
Nill and Szlach\'anyi in \cite{bohm1}, the notion of crossed product
can be adapted to the weak setting. In the Hopf algebra world,
crossed products appear as a generalization of semi-direct products
of groups to the context of Hopf algebras \cite{Molnar,bcm}, and are
closely connected   with  cleft extensions and Galois extensions of
Hopf algebras \cite{blat-susan, doi3}. In \cite{tb-crpr}
Brzezi\'nski gave an interesting approach that generalizes several
types of crossed products, even the ones given for braided Hopf
algebras by Majid \cite{maj2} and Guccione and Guccione in
\cite{Guccione}. On the other hand, in \cite{LS} we can find a
general and categorical theory, the theory of wreath products, that
contains as a particular instance the crossed structures presented
by Brzezi\'nski.

The key  to extend the crossed product constructions presented in
the previous paragraph  to the weak setting is the  use of
idempotent morphisms combined with the ideas in \cite{tb-crpr}. In
\cite{nmra4} the authors defined a product on $A\otimes V$, for an
algebra $A$ and an object $V$ both living in a strict monoidal
category $\mathcal C$ where every idempotent splits. In order to
obtain that product we must consider two morphisms
$\psi_{V}^{A}:V\otimes A\rightarrow A\otimes V$ and
$\sigma_{V}^{A}:V\otimes V\rightarrow A\otimes V$ that satisfy some
twisted-like and cocycle-like conditions. Associated to these
morphisms it is possible to define an idempotent morphism
$\nabla_{A\otimes V}:A\otimes V\rightarrow A\otimes V$ and the image
of $\nabla_{A\otimes V}$ inherits the associative product from
$A\otimes V$. In order to define a unit for $Im (\nabla_{A\otimes
V})$, and hence to obtain an algebra structure, we require the
existence of  a preunit $\nu:K\rightarrow A\otimes V$.  In
\cite{mra-preunit} we can find a characterization  of weak crossed
products with a preunit as associative products on $A\otimes V$ that
are morphisms of left $A$-modules with preunit. Finally, it is
convenient to observe that, if the preunit is an unit, the
idempotent becomes the identity and we recover the classical
examples of the Hopf algebra setting. The theory presented in
\cite{nmra4, mra-preunit} contains as a particular instance the one
developed by Brzezi\'nski in \cite{tb-crpr}. There are many other
examples of this theory like the weak smash product given by
Caenepeel and De Groot in \cite{caengroot}, the theory of wreath
products presented in \cite{LS} and the weak crossed products for
weak bialgebras given in \cite{ana1}. Recently, G. B\"ohm showed in
\cite{bohm} that a monad in the weak version of the Lack and
Street's 2-category of monads in a 2-category is identical to a
crossed product system in the sense of \cite{nmra4} and also in
\cite{mra-partial-unif} we can find that unified crossed products
\cite{AM1} and partial crossed products \cite{partial} are
particular instances of weak crossed products.

Then, if in the Hopf algebra setting the second cohomology group
classifies crossed products of $H$ with a commutative left
$H$-module algebra $A$, what about the weak setting? The answer to
this question is the main motivation of this paper. More precisely,
we show that if $H$ is a cocommutative weak Hopf algebra and $A$ is
a commutative left $H$-module algebra, all the weak crossed products
defined in $A\ot H$ with a common preunit can be described by the
second cohomology group of a new cohomology that we call the
Sweedler cohomology of a weak Hopf algebra with coefficients in $A$.

The paper is organized  as follows: In Section 1 after recalling the
basic properties of weak Hopf algebras, we introduce the notion of
weak $H$-module algebra and  define the cosimplicial complex
$Reg_{\varphi_{A}}(H^{\bullet}, A)$ for a cocommutative weak Hopf
algebra $H$ and a commutative left $H$-module algebra $A$. Then, we
introduce the Sweedler cohomology of $H$ with coefficients in $A$ as
the one defined by the associated cochain complex.  In the second
section we present the results about the characterization of weak
crossed products induced by morphisms $\sigma\in
Reg_{\varphi_{A}}(H^{2}, A)$ proving that the twisted condition, and
the cocycle condition of the general theory of weak crossed products
can be reduced to twisted an 2-cocycle for the action and $\sigma$.
Also, in this section we introduce the normal condition that permits
to obtain a preunit in the weak crossed product induced by the
morphism $\sigma$. Finally, in the third section we characterize the
equivalence between two weak crossed products obtaining the main
result of this paper that assures the following: There is a
bijective correspondence between $H^{2}_{\varphi_{A}}(H,A)$ and the
equivalence classes of weak crossed products of $A\ot_{\alpha} H$
where $\alpha: H\ot H\rightarrow A$ satisfies the 2-cocycle
 and the normal conditions.

\section{The  Sweedler cohomology in a weak setting}

From now on ${\mathcal C}$ denotes a strict symmetric category with
tensor product denoted by $\ot$ and unit object $K$. With $c$ we
will denote the natural isomorphism of symmetry and we also assume
that ${\mathcal C}$ has equalizers. Then, under these conditions,
every idempotent morphism $q:Y\rightarrow Y$ splits, i.e., there
exist an object $Z$ and morphisms $i:Z\rightarrow Y$ and
$p:Y\rightarrow Z$ such that $q=i\circ p$ and $p\circ i =id_{Z}$. We
denote the class of objects of ${\mathcal C}$ by $\vert {\mathcal C}
\vert $ and for each object $M\in \vert {\mathcal C}\vert$, the
identity morphism by $id_{M}:M\rightarrow M$. For simplicity of
notation, given objects $M$, $N$, $P$ in ${\mathcal C}$ and a
morphism $f:M\rightarrow N$, we write $P\ot f$ for $id_{P}\ot f$ and
$f \ot P$ for $f\ot id_{P}$.

An algebra in ${\mathcal C}$ is a triple $A=(A, \eta_{A}, \mu_{A})$
where $A$ is an object in ${\mathcal C}$ and
 $\eta_{A}:K\rightarrow A$ (unit), $\mu_{A}:A\otimes A
\rightarrow A$ (product) are morphisms in ${\mathcal C}$ such that
$\mu_{A}\circ (A\otimes \eta_{A})=id_{A}=\mu_{A}\circ
(\eta_{A}\otimes A)$, $\mu_{A}\circ (A\otimes \mu_{A})=\mu_{A}\circ
(\mu_{A}\otimes A)$. We will say that an algebra $A$ is commutative
if $\mu_{A}\co c_{A,A}=\mu_{A}$.

Given two algebras $A= (A, \eta_{A}, \mu_{A})$ and $B=(B, \eta_{B},
\mu_{B})$, $f:A\rightarrow B$ is an algebra morphism if
$\mu_{B}\circ (f\otimes f)=f\circ \mu_{A}$ and  $ f\circ \eta_{A}=
\eta_{B}.$

If $A$, $B$ are algebras in ${\mathcal C}$, the object $A\otimes B$
is an algebra in
 ${\mathcal C}$ where
$\eta_{A\otimes B}=\eta_{A}\otimes \eta_{B}$ and
\begin{equation}
\label{it-prod} \mu_{A\otimes B}=(\mu_{A}\otimes \mu_{B})\circ
(A\otimes c_{B,A}\otimes B).
\end{equation}

A coalgebra in ${\mathcal C}$ is a triple ${D} = (D,
\varepsilon_{D}, \delta_{D})$ where $D$ is an object in ${\mathcal
C}$ and $\varepsilon_{D}: D\rightarrow K$ (counit),
$\delta_{D}:D\rightarrow D\otimes D$ (coproduct) are morphisms in
${\mathcal C}$ such that $(\varepsilon_{D}\otimes D)\circ
\delta_{D}= id_{D}=(D\otimes \varepsilon_{D})\circ \delta_{D}$,
$(\delta_{D}\otimes D)\circ \delta_{D}=
 (D\otimes \delta_{D})\circ \delta_{D}.$ We will say that $D$ is
 cocommutative if $c_{D,D}\co \delta_{D}=\delta_{D}$ holds.

 If ${D} = (D, \varepsilon_{D},
 \delta_{D})$ and
${ E} = (E, \varepsilon_{E}, \delta_{E})$ are coalgebras,
$f:D\rightarrow E$ is a coalgebra morphism if $(f\otimes f)\circ
\delta_{D} =\delta_{E}\circ f$ and $\varepsilon_{E}\circ f
=\varepsilon_{D}.$

When $D$, $E$ are coalgebras in ${\mathcal C}$, $D\otimes E$ is a
coalgebra in ${\mathcal C}$ where $\varepsilon_{D\otimes
E}=\varepsilon_{D}\otimes \varepsilon_{E}$ and
\begin{equation}
\label{it-coprod} \delta_{D\otimes E}=(D\otimes c_{D,E}\otimes
E)\circ( \delta_{D}\otimes \delta_{E}).
\end{equation}

If $A$ is an algebra, $B$ is a coalgebra and $\alpha:B\rightarrow
A$, $\beta:B\rightarrow A$ are morphisms, we define the convolution
product by
$$\alpha\wedge \beta=\mu_{A}\circ (\alpha\otimes
\beta)\circ \delta_{B}.$$

Let  $A$ be an algebra. The pair $(M,\phi_{M})$ is a right
$A$-module if $M$ is an object in ${\mathcal C}$ and
$\phi_{M}:M\otimes A\rightarrow M$ is a morphism in ${\mathcal C}$
satisfying $\phi_{M}\circ(M\otimes \eta_{A})=id_{M}$, $\phi_{M}\circ
(\phi_{M}\otimes A)=\phi_{M}\circ (M\otimes \mu_{A})$. Given two
right ${A}$-modules $(M,\phi_{M})$ and $(N,\phi_{N})$,
$f:M\rightarrow N$ is a morphism of right ${A}$-modules if
$\phi_{N}\circ (f\otimes A)=f\circ \phi_{M}$. In a similar way we
can define the notions of left $A$-module and morphism of left
$A$-modules. In this case we denote the left action by
$\varphi_{M}$.

 Let  $C$ be a coalgebra. The pair
$(M,\rho_{M})$ is a right $C$-comodule if $M$ is an object in
${\mathcal C}$ and $\rho_{M}:M\rightarrow M\ot C$ is a morphism in
${\mathcal C}$ satisfying $(M\otimes \varepsilon_{C})\co
\rho_{M}=id_{M}$, $(M\ot \rho_{M})\co \rho_{M}=(M\ot \delta_{C})\co
\rho_{M}$. Given two right ${C}$-comodules $(M,\rho_{M})$ and
$(N,\rho_{N})$, $f:M\rightarrow N$ is a morphism of right
${C}$-comodules if $(f\otimes C)\co \rho_{M}=\rho_{N}\co f$. In a
similar way we can define the notions of left $C$-comodule and
morphism of left $C$-comodules. In this case we denote the left
action by $\varrho_{M}$.

By weak Hopf algebras  we understand the objects introduced in
\cite{bohm1}, as a generalization of ordinary Hopf algebras. Here we
recall the definition of these objects in the symmetric monoidal
setting.

\begin{defin}
\label{wha} {\rm A weak Hopf algebra $H$  is an object in ${\mathcal
C}$ with an algebra structure $(H, \eta_{H},\mu_{H})$ and a
coalgebra structure $(H, \varepsilon_{H},\delta_{H})$ such that the
following axioms hold:
\begin{itemize}

\item[(a1)] $\delta_{H}\circ \mu_{H}=(\mu_{H}\otimes \mu_{H})\circ
\delta_{H\otimes H},$

\item[(a2)]$\varepsilon_{H}\circ \mu_{H}\circ
(\mu_{H}\otimes H)=(\varepsilon_{H}\otimes \varepsilon_{H})\circ
(\mu_{H}\otimes \mu_{H})\circ (H\otimes \delta_{H}\otimes H)$
\item[ ]$=(\varepsilon_{H}\otimes \varepsilon_{H})\circ
(\mu_{H}\otimes \mu_{H})\circ (H\otimes
(c_{H,H}\circ\delta_{H})\otimes H),$

\item[(a3)]$(\delta_{H}\otimes H)\circ \delta_{H}\circ
\eta_{H}=(H\otimes \mu_{H}\otimes H)\circ (\delta_{H}\otimes
\delta_{H})\circ (\eta_{H}\otimes \eta_{H})$ \item[ ]$=(H\otimes
(\mu_{H}\circ c_{H,H})\otimes H)\circ (\delta_{H}\otimes
\delta_{H})\circ (\eta_{H}\otimes \eta_{H}).$

\item[(a4)] There exists a morphism $\lambda_{H}:H\rightarrow H$
in ${\mathcal C}$ (called the antipode of $H$) satisfying:

\begin{itemize}

\item[(a4-1)] $id_{H}\wedge \lambda_{H}=((\varepsilon_{H}\circ
\mu_{H})\otimes H)\circ (H\otimes c_{H,H})\circ ((\delta_{H}\circ
\eta_{H})\otimes H),$

\item[(a4-2)] $\lambda_{H}\wedge
id_{H}=(H\otimes(\varepsilon_{H}\circ \mu_{H}))\circ (c_{H,H}\otimes
H)\circ (H\otimes (\delta_{H}\circ \eta_{H})),$

\item[(a4-3)]$\lambda_{H}\wedge id_{H}\wedge
\lambda_{H}=\lambda_{H}.$
\end{itemize}

\end{itemize}
}
\end{defin}

\begin{apart}
{\rm If $H$ is a weak Hopf algebra in ${\mathcal C}$, the antipode
$\lambda_{H}$ is unique, antimultiplicative, anticomultiplicative
and leaves the unit  and the counit invariant:
\begin{equation}
\label{anti} \lambda_{H}\circ \mu_{H}=\mu_{H}\circ
(\lambda_{H}\otimes \lambda_{H})\circ
c_{H,H};\;\;\;\;\delta_{H}\circ \lambda_{H}=c_{H,H}\circ
(\lambda_{H}\otimes \lambda_{H})\circ \delta_{H};
\end{equation}
\begin{equation}
\label{invariant-unit-counit} \lambda_{H}\circ
\eta_{H}=\eta_{H};\;\;\;\;\varepsilon_{H}\circ
\lambda_{H}=\varepsilon_{H}.
\end{equation}

If we define the morphisms $\Pi_{H}^{L}$ (target), $\Pi_{H}^{R}$
(source), $\overline{\Pi}_{H}^{L}$ and $\overline{\Pi}_{H}^{R}$ by
\begin{itemize}

\item[ ]$\Pi_{H}^{L}=((\varepsilon_{H}\circ \mu_{H})\otimes
H)\circ (H\otimes c_{H,H})\circ ((\delta_{H}\circ \eta_{H})\otimes
H),$

\item[ ]$\Pi_{H}^{R}=(H\otimes(\varepsilon_{H}\circ
\mu_{H}))\circ (c_{H,H}\otimes H)\circ (H\otimes (\delta_{H}\circ
\eta_{H})),$

\item[ ]$\overline{\Pi}_{H}^{L}=(H\otimes
(\varepsilon_{H}\circ \mu_{H}))\circ ((\delta_{H}\circ
\eta_{H})\otimes H),$

\item[ ]$\overline{\Pi}_{H}^{R}=((\varepsilon_{H}\circ \mu_{H})\otimes
H)\circ(H\otimes (\delta_{H}\circ \eta_{H})),$

\end{itemize}
it is straightforward to show (see \cite{bohm1}) that they are
idempotent and $\Pi_{H}^{L}$, $\Pi_{H}^{R}$ satisfy the equalities
\begin{equation}
\Pi_{H}^{L}=id_{H}\wedge
\lambda_{H};\;\;\;\;\;\Pi_{H}^{R}=\lambda_{H}\wedge id_{H}.
\end{equation}
and then
\begin{equation}
\label{id} \Pi_{H}^{L}\wedge \Pi_{H}^{L}=\Pi_{H}^{L},\;\;\;
\Pi_{H}^{R}\wedge \Pi_{H}^{R}=\Pi_{H}^{R}.
\end{equation}

Moreover, we have that
\begin{equation}
\label{lastone} \Pi_{H}^{L}\circ
\overline{\Pi}_{H}^{L}=\Pi_{H}^{L};\;\;\;\; \Pi_{H}^{L}\circ
\overline{\Pi}_{H}^{R}=\overline{\Pi}_{H}^{R};\;\;\;\;
\Pi_{H}^{R}\circ
\overline{\Pi}_{H}^{L}=\overline{\Pi}_{H}^{L};\;\;\;\;
\Pi_{H}^{R}\circ \overline{\Pi}_{H}^{R}=\Pi_{H}^{R};
\end{equation}
\begin{equation}
\overline{\Pi}_{H}^{L}\circ
\Pi_{H}^{L}=\overline{\Pi}_{H}^{L};\;\;\;\;
\overline{\Pi}_{H}^{L}\circ \Pi_{H}^{R}=\Pi_{H}^{R};\;\;\;\;
\overline{\Pi}_{H}^{R}\circ \Pi_{H}^{L}=\Pi_{H}^{L};\;\;\;\;
\overline{\Pi}_{H}^{R}\circ \Pi_{H}^{R}=\overline{\Pi}_{H}^{R}.
\end{equation}

For the morphisms target an source we have the following identities:
\begin{equation}
\label{deltamu} \Pi^{L}_{H}\circ \mu_{H}\circ (H\ot
\Pi^{L}_{H})=\Pi^{L}_{H}\circ \mu_{H},\;\;\; \Pi^{R}_{H}\circ
\mu_{H}\circ (\Pi^{R}_{H}\ot H)=\Pi^{R}_{H}\circ \mu_{H},
\end{equation}
\begin{equation}
\label{deltapi} (H\ot \Pi^{L}_{H})\circ \delta_{H}\circ
\Pi^{L}_{H}=\delta_{H}\circ \Pi^{L}_{H},\;\;\;( \Pi^{R}_{H}\ot
H)\circ \delta_{H}\circ \Pi^{R}_{H}=\delta_{H}\circ \Pi^{R}_{H},
\end{equation}
\begin{equation}
\label{delta-pi1} \mu_{H}\co (H\ot \Pi_{H}^{L})=((\varepsilon_{H}\co
\mu_{H})\ot H)\co (H\ot c_{H,H})\co (\delta_{H}\ot H),
\end{equation}
\begin{equation}
\label{delta-pi2} (H\ot \Pi^{L}_{H})\co \delta_{H}= (\mu_{H}\ot
H)\co (H\ot c_{H,H})\co ((\delta_{H}\co \eta_{H})\ot H),
\end{equation}
\begin{equation}
\label{delta-pi3} \mu_{H}\co (\Pi_{H}^{R} \ot H)=(H\ot
(\varepsilon_{H}\co \mu_{H}))\co (c_{H,H}\ot H)\co (H\ot \delta_{H})
\end{equation}
\begin{equation}
\label{delta-pi4} (\Pi_{H}^{R} \ot H)\co \delta_{H}=(H\ot
\mu_{H})\co (c_{H,H}\ot H)\co (H\ot (\delta_{H}\co \eta_{H}))
\end{equation}
and
\begin{equation}
\label{delta-pi11} \mu_{H}\co (\overline{\Pi}_{H}^{R}\ot
H)=((\varepsilon_{H}\co \mu_{H})\ot H)\co (H\ot \delta_{H}),
\end{equation}
\begin{equation}
\label{delta-pi21} \mu_{H}\co (H\ot \overline{\Pi}_{H}^{L})=(H\ot
(\varepsilon_{H}\co \mu_{H}))\co (\delta_{H}\ot H),
\end{equation}
\begin{equation}
\label{delta-pi31} (\overline{\Pi}_{H}^{L}\ot H)\co \delta_{H}=(H\ot
\mu_{H})\co ((\delta_{H}\co \eta_{H})\ot H),
\end{equation}
\begin{equation}
\label{delta-pi41} (H\ot \overline{\Pi}_{H}^{R})\co
\delta_{H}=(\mu_{H}\ot H)\co (H\ot (\delta_{H}\co \eta_{H})),
\end{equation}

Finally, if $H$ is (co)commutative we have that $\lambda_{H}$ is an
isomorphism and $\lambda_{H}^{-1}=\lambda_{H}$.

}

\end{apart}

\begin{ej}
\label{q-group} {\rm As group algebras and their duals are the
natural examples of Hopf algebras, groupoid algebras and their duals
provide examples of weak Hopf algebras. Recall that a groupoid $G$
is simply a category in which every morphism is an isomorphism. In
this example, we consider finite groupoids, i.e. groupoids with a
finite number of objects. The set of objects of $G$ will be denoted
by $G_{0}$ and the set of morphisms by $G_{1}$. The identity
morphism on $x\in G_{0}$ will also be denoted by $id_{x}$ and for a
morphism $\sigma:x\rightarrow y$ in $G_{1}$, we write $s(\sigma)$
and $t(\sigma)$, respectively for the source and the target of
$\sigma$.

Let $G$ be a groupoid, and $R$ a commutative ring. The groupoid
algebra is the direct product
$$RG=\bigoplus_{\sigma\in G_{1}}R\sigma$$

with the product of two morphisms being equal to their composition
if the latter is defined and $0$  otherwise, i.e.
$\sigma\tau=\sigma\circ \tau$ if $s(\sigma)=t(\tau)$ and
$\sigma\tau=0$ if $s(\sigma)\neq t(\tau)$. The unit element is
$1_{RG}=\sum_{x\in G_{0}}id_{x}$. Then $RG$ is a cocommutative weak
Hopf algebra, with coproduct $\delta_{RG}$, counit
$\varepsilon_{RG}$ and antipode $\lambda_{RG}$ given by the formulas
$\delta_{RG}(\sigma)=\sigma\otimes \sigma,$ $\varepsilon_{RG}
(\sigma)=1,$ and $\lambda_{RG}(\sigma)=\sigma^{-\dot{}1}.$ For the
weak Hopf algebra $RG$ the morphisms target and source are
respectively,
$$\Pi_{RG}^{L}(\sigma)=id_{t(\sigma)},\;\;\;
\Pi_{RG}^{R}(\sigma)=id_{s(\sigma)}.$$}
\end{ej}

\begin{defin}
\label{weak-H-mod} {\rm Let $H$ be a weak Hopf algebra. We will say
that $A$ is a weak left $H$-module algebra if there exists a
morphism $\varphi_{A}:H\ot A\rightarrow A$ satisfying:

\begin{itemize}

\item[(b1)] $\varphi_{A}\co (\eta_{H}\ot A)=id_{A}$.

\item[(b2)] $\varphi_{A}\co (H\ot \mu_{A})=\mu_{A}\co
(\varphi_{A}\ot \varphi_{A})\co (H\ot c_{H,A}\ot A)\co
(\delta_{H}\ot A\ot A).$

\item[(b3)] $\varphi_{A}\co (\mu_{H}\ot \eta_{A})=\varphi_{A}\co (H\ot
(\varphi_{A}\co (H\ot \eta_{A}))).$

\end{itemize}

and any of the following equivalent conditions holds:

\begin{itemize}

\item[(b4)] $\varphi_{A}\co (\Pi^{L}_{H}\ot A)=\mu_{A}\co ((\varphi_{A}\co (H\ot
\eta_{A})\ot A).$

\item[(b5)] $\varphi_{A}\co (\overline{\Pi}^{L}_{H}\ot A)=\mu_{A}\co c_{A,A}\co
((\varphi_{A}\co (H\ot \eta_{A})\ot A).$

\item[(b6)] $\varphi_{A}\co (\Pi^{L}_{H}\ot \eta_{A})=\varphi_{A}\co (H\ot \eta_{A}).$

\item[(b7)] $\varphi_{A}\co (\overline{\Pi}^{L}_{H}\ot\eta_{A})=
\varphi_{A}\co (H\ot \eta_{A}).$

\item[(b8)] $\varphi_{A}\co (H\ot (\varphi_{A}\co (H\ot
\eta_{A})))=((\varphi_{A}\co (H\ot \eta_{A}))\ot (\varepsilon_{H}\co
\mu_{H}))\co (\delta_{H}\ot H).$

\item[(b9)] $\varphi_{A}\co (H\ot (\varphi_{A}\co (H\ot
\eta_{A})))=((\varepsilon_{H}\co \mu_{H})\ot (\varphi_{A}\co (H\ot
\eta_{A})))\co (H\ot c_{H,H})\co (\delta_{H}\ot H).$

\end{itemize}

If we replace (b3) by
\begin{itemize}

\item[(b3-1)] $\varphi_{A}\co (\mu_{H}\ot A)=\varphi_{A}\co (H\ot
\varphi_{A})$

\end{itemize}
we will say that $(A, \varphi_{A})$ is a left $H$-module algebra.

}
\end{defin}

\begin{notac}
\label{notation-1} {\rm Let $H$ be a weak Hopf algebra.  For $n\geq
1$, we  denote by  $H^n$ the $n$-fold tensor power $H\ot \cdots \ot
H$.  By $H^{0}$ we denote the unit object of ${\mathcal C}$, i.e.
$H^{0}=K$.

If $n\geq 2$,  $m_{H}^{n}$ denotes the morphism
$$m_{H}^{n}:H^n\rightarrow H$$
defined by $m_{H}^{2}=\mu_{H}$ and by
$$m_{H}^{3}=m_{H}^{2}\co (H\ot \mu_{H}),\cdots, m_{H}^{n}=m_{H}^{n-1}\co (H^{n-2}\ot \mu_{H})
$$
for $k> 2$. Note that by the associativity of $\mu_{H}$ we have
$$m_{H}^{n}=m_{H}^{n-1}\co (\mu_{H}\ot H^{n-2}).
$$

Let $(A,\varphi_{A})$ be a weak left $H$-module algebra and $n\geq
1$. With  $\varphi^{n}_{A}$ we will denote the morphism
$$\varphi^{n}_{A}:H^n\ot A\rightarrow A$$
defined as $\varphi^{1}_{A}=\varphi_{A}$ and
$\varphi^{n}_{A}=\varphi_{A}\co (H\ot \varphi^{n-1}_{A})$. If $n>
1$, we have that
\begin{equation}
\label{varphi-eta} \varphi_{A}\co (m_{H}^{n}\ot
\eta_{A})=\varphi_{A}^{n-1}\co (H^{n-1}\ot (\varphi_{A}\co
(H\ot\eta_{A}))
\end{equation}
holds. In what follows,  we denote the morphism $\varphi_{A}\co
(m_{H}^{n}\ot \eta_{A})$ by $u_{n}$ and the morphism $\varphi_{A}\co
(H\ot \eta_{A})$ by $u_{1}$. Note that, by (b3) of Definition
\ref{weak-H-mod},  for $n\geq 2$,
\begin{equation}
\label{varphi-eta-1} u_{n}=\varphi_{A}^{n-1}\co (H^{n-1}\ot u_{1}).
\end{equation}

Finally, with $\delta_{H^{n}}$ we denote the coproduct defined in
(\ref{it-coprod}) for the coalgebra $H^{n}$. Then,
\begin{equation}
\label{it-coprod-H} \delta_{H^{n}}=\delta_{H^{k}\ot
H^{n-k}}=\delta_{H^{n-k}\ot H^{k}},
\end{equation}
for $k\in\{1,\dots,n-1\}$.
}
\end{notac}

\begin{prop}
\label{cocommutative-properties} Let $H$ be a cocommutative weak
Hopf algebra. The following identities hold.
\begin{itemize}

\item[(i)] $\delta_{H}\co \Pi_{H}^{I}=(\Pi_{H}^{I}\ot
\Pi_{H}^{I})\co \delta_{H}$ for $I\in\{L,R\}.$

\item[(ii)] $ (\Pi_{H}^{I}\ot H)\co \delta_{H}\co \Pi_{H}^{J}=
(H\ot \Pi_{H}^{I})\co \delta_{H}\co \Pi_{H}^{J}=\delta_{H}\co
\Pi_{H}^{J}$, for $I,J\in\{L,R\}.$

\item[(iii)] $(\Pi_{H}^{L}\ot H)\co \delta_{H}\co \mu_{H}=(\Pi_{H}^{L}\ot \mu_{H})\co
(\delta_{H}\ot H) $.

\item[(iv)] $(H\ot \Pi_{H}^{R})\co \delta_{H}\co \mu_{H}=( \mu_{H}\ot \Pi_{H}^{R})\co
(H\ot \delta_{H})$.

\end{itemize}

\end{prop}

{\em Proof}: First note that if $H$ is cocommutative
$\Pi_{H}^{I}=\overline{\Pi}_{H}^{I}$ for $I\in\{L,R\}.$ The proof
for (i) with $I=L$ follows by

\begin{itemize}

\item[ ]$\hspace{0.38cm}\delta_{H}\co \Pi_{H}^{L} $

\item [ ]$= \mu_{H\ot H}\co (\delta_{H}\ot (\delta_{H}\co \lambda_{H}))\co \delta_{H}$

\item [ ]$= \mu_{H\ot H}\co (\delta_{H}\ot (c_{H,H}\co (\lambda_{H}\ot \lambda_{H})\co \delta_{H}))\co \delta_{H}  $

\item [ ]$=(\mu_{H}\ot \Pi_{H}^{L} )\co (H\ot c_{H,H})\co (c_{H,H}\ot \lambda_{H})
\co (H\ot \delta_{H})\co \delta_{H}  $

\item [ ]$=(\Pi_{H}^{L}\ot
\Pi_{H}^{L})\co \delta_{H} $

\end{itemize}

where the first equality follows by (a1) of Definition \ref{wha},
the second by  the antimultiplicative property of $\lambda_{H}$, the
third one relies on the naturality of $c$, the coassociativity of
$\delta_{H}$ and the cocommutativity of $H$. Finally, the last one
follows by the cocommutativity of $H$ and the naturality of $c$.

The proof for $I=R$ is similar.

Note that, by (i) and the idempotent property of $\Pi_{H}^{I}$, we
have (ii) for $I=J$. If $I=L$ and $J=R$, by (\ref{lastone}),  we
have
$$(\Pi_{H}^{L}\ot H)\co \delta_{H}\co \Pi_{H}^{R}=((\Pi_{H}^{L}\co
\Pi_{H}^{R})\ot \Pi_{H}^{R})\co \delta_{H}=((\Pi_{H}^{L}\co
\overline{\Pi}_{H}^{R})\ot \Pi_{H}^{R})\co \delta_{H}$$
$$=(\overline{\Pi}_{H}^{R}\ot \Pi_{H}^{R})\co \delta_{H}=(\Pi_{H}^{R}\ot \Pi_{H}^{R})\co \delta_{H}=\delta_{H}\co
\Pi_{H}^{R}.$$

The proof for  $I=R$ and $J=L$ is similar. On the other hand, by the
usual arguments, we get (iii):
$$(\Pi_{H}^{L}\ot H)\co \delta_{H}\co \mu_{H} = (\overline{\Pi}_{H}^{L}\ot H)\co \delta_{H}\co \mu_{H}
= (H\ot \mu_{H})\co ((\delta_{H}\co \eta_{H})\ot \mu_{H})  $$
$$=(\overline{\Pi}_{H}^{L}\ot \mu_{H})\co (\delta_{H}\ot H)
=(\Pi_{H}^{L}\ot \mu_{H})\co (\delta_{H}\ot H). $$

The proof of the equality (iv) follows a similar pattern and we
leave the details to the reader.

\begin{prop}
\label{cocommutative-properties-2} Let $H$ be a cocommutative weak
Hopf algebra. The following identities hold.
\begin{itemize}

\item[(i)] $ \delta_{H^{2}}\co \delta_{H}=(\delta_{H}\ot \delta_{H})\co \delta_{H} .$

\item[(ii)] $\delta_{H^{n+1}}\co ( H^{i}\ot \delta_{H}\ot H^{n-i-1})=(H^{i}\ot
\delta_{H}\ot H^{n-1}\ot \delta_{H}\ot H^{n-i-1}) \co
\delta_{H^{n}}$ for  $n\geq 2$ and $i\in\{0,\cdots,n-1\}.$

\item[(iii)] $\delta_{H^{n}}\co (H^{i}\ot \Pi_{H}^{I}\ot H^{n-i-1}) =
(H^{i}\ot \Pi_{H}^{I}\ot H^{n-1}\ot \Pi_{H}^{I}\ot H^{n-i-1})\co
\delta_{H^{n}} $ for $I\in\{L,R\}$, $n\geq 2$ and
$i\in\{0,\cdots,n-1\}$.

\item[(iv)] $ \delta_{H^{n+1}}\co ( H^{i}\ot ((\Pi_{H}^{I}\ot H)\co \delta_{H})\ot H^{n-i-1})$
$$=(H^{i}\ot ((\Pi_{H}^{I}\ot H)\co \delta_{H})\ot H^{n-1}\ot
((\Pi_{H}^{I}\ot H)\co \delta_{H})\ot H^{n-i-1}) \co
\delta_{H^{n}}$$ for $I\in\{L,R\}$, $n\geq 2$ and
$i\in\{0,\cdots,n-1\}$.

\item[(v)] $ \delta_{H^{n+1}}\co ( H^{i}\ot ((H\ot \Pi_{H}^{I})\co \delta_{H})\ot H^{n-i-1})$
$$=(H^{i}\ot ((H\ot \Pi_{H}^{I})\co \delta_{H})\ot H^{n-1}\ot
((H\ot \Pi_{H}^{I})\co \delta_{H})\ot H^{n-i-1}) \co
\delta_{H^{n}}$$ for $I\in\{L,R\}$, $n\geq 2$ and
$i\in\{0,\cdots,n-1\}$.

\end{itemize}

\end{prop}

{\em Proof}: The assertion (i) follows by the coassociativity of
$\delta_{H}$ and the cocommutativity of $H$. Indeed:
$$\delta_{H^{2}}\co \delta_{H} = (H\ot (c_{H,H}\co \delta_{H})\ot
H)\co (\delta_{H}\ot H)\co \delta_{H}=(\delta_{H}\ot \delta_{H})\co
\delta_{H}.$$

The proof for (ii)  can be obtained using (i) and mathematical
induction. Also, by this method and Proposition
\ref{cocommutative-properties} we obtain (iii), (iv) and (iv).

\begin{rem}
\label{coalgebra-structure} {\rm If $H$ is a weak Hopf algebra, we
denote by $H_{L}$ the object such that $p_{L}\co i_{L}=id_{H_{L}}$
where $i_{L}$, $p_{L}$ are the injection and the projection
associated to the target morphism $\Pi_{H}^{L}$. If $H$ is
cocommutative, by (i) of Proposition \ref{cocommutative-properties},
we have that $H_{L}$ is a coalgebra and the morphisms $i_{L}$,
$p_{L}$ are coalgebra morphisms for $\delta_{H_{L}}=(p_{L}\ot
p_{L})\co \delta_{H}\co i_{L}$ and
$\varepsilon_{H_{L}}=\varepsilon_{H}\co i_{L}$. Therefore,
$\delta_{H_{L}}\co p_{L}=(p_{L}\ot p_{L})\co \delta_{H}$ and
$\varepsilon_{H_{L}}\co p_{L}=\varepsilon_{H}.$

}
\end{rem}

\begin{prop}
\label{delta-mu-it} Let $H$ be a weak Hopf algebra. Then, if $n\geq
3$ the following equality holds.
\begin{equation}
\label{prod-delta}
(H^{i-1}\ot \mu_{H}\ot H^{n-i-1}\ot H^{i-1}\ot \mu_{H}\ot
H^{n-i-1})\co \delta_{H^{n}}=\delta_{H^{n-1}}\co (H^{i-1}\ot
\mu_{H}\ot H^{n-i-1}),
\end{equation}
for all $i\in\{1,\cdots, n-1\}.$
\end{prop}

{\em Proof}: First note that, by (a1) of Definition \ref{wha}, we
have $(\mu_{H}\ot \mu_{H})\co \delta_{H^{2}}=\delta_{H}\co \mu_{H}.$
Then, using this identity we have:
$$(\mu_{H}\ot
H^{n-2}\ot \mu_{H}\ot H^{n-2})\co \delta_{H^{n}}=(\mu_{H}\ot
H^{n-2}\ot \mu_{H}\ot H^{n-2})\co (H^2\ot c_{H^2,H^{n-2}}\ot
H^{n-2})\co (\delta_{H^{2}}\ot \delta_{H^{n-2}})=$$
$$(H\ot c_{H,H^{n-2}}\ot H^{n-2})\co (((\mu_{H}\ot \mu_{H})\co
\delta_{H^{2}})\ot \delta_{H^{n-2}})=\delta_{H^{n-1}}\co (\mu_{H}\ot
H^{n-2}).$$

Then, as a consequence, we have

\begin{itemize}
\item[ ]$\hspace{0.38cm}(H^{i-1}\ot \mu_{H}\ot H^{n-i-1}\ot H^{i-1}\ot \mu_{H}\ot
H^{n-i-1})\co \delta_{H^{n}} $

\item [ ]$= (H^{i-1}\ot \mu_{H}\ot
H^{n-i-1}\ot H^{i-1}\ot \mu_{H}\ot H^{n-i-1})\co  $
\item[ ]$\hspace{0.38cm} (H^{i-1}\ot H^2\ot
c_{H^{i-1}, H^{n-i-1}}\ot H^2\ot H^{n-i-1})\co (H^{i-1}\ot
c_{H^{i-1}, H^2}\ot H^{n-i-1}\ot H^2\ot H^{n-i-1})\co$
\item[ ]$\hspace{0.38cm} (\delta_{H^{i-1}}\ot \delta_{H^{n-i+1}} )$

\item [ ]$=  (H^{i-1}\ot H\ot
c_{H^{i-1}, H^{n-i-1}}\ot H\ot H^{n-i-1})\co (H^{i-1}\ot c_{H^{i-1},
H}\ot H^{n-i-1}\ot H\ot H^{n-i-1})\co$
\item[ ]$\hspace{0.38cm}
(\delta_{H^{i-1}}\ot ((\mu_{H}\ot H^{n-i-1}\ot \mu_{H}\ot
H^{n-i-1})\co \delta_{H^{n-i+1}}) )$

\item [ ]$=  (H^{i-1}\ot
c_{H^{i-1}, H^{n-i}}\ot H^{n-i})\co (\delta_{H^{i-1}}\ot (
\delta_{H^{n-i}}\co (\mu_{H}\ot H^{n-i-1})))$

\item [ ]$=\delta_{H^{n-1}}\co (H^{i-1}\ot
\mu_{H}\ot H^{n-i-1}).$

\end{itemize}

\begin{prop}
\label{iterations2} Let $H$ be a  weak Hopf algebra. The following
identity holds for $n\geq 2$.
\begin{equation}
\label{m-delta} \delta_{H}\co m^{n}_{H}=(m^{n}_{H}\ot
m^{n}_{H})\co\delta_{H^{n}}.
\end{equation}

\end{prop}

{\em Proof}: As in the previous proposition we proceed by induction.
Obviously the equality (\ref{m-delta}) holds  for $n=2$. If we
assume that it is true for $n=k$, it is true for $n=k+1$ because:
$$(m^{k+1}_{H}\ot
m^{k+1}_{H})\co\delta_{H^{k+1}}=((\mu_{H}\co (m^{k}_{H}\ot H))\ot
(\mu_{H}\co (m^{k}_{H}\ot H)))\co \delta_{H^{k}\ot H}$$
$$=\mu_{H\ot
H}\co (((m^{k}_{H}\ot m^{k}_{H})\co \delta_{H^{k}}))\ot
\delta_{H})=\delta_{H}\co \mu_{H}\co (m^{k}_{H}\ot H)= \delta_{H}\co
m^{k+1}_{H}.$$

\begin{prop}
\label{unit} Let $H$ be a weak Hopf algebra and $(A,\varphi_{A})$ be
a weak left $H$-module algebra. Then, if $n\geq 1$, the equality
\begin{equation}
\label{neutro} u_{n}\wedge u_{n}=u_{n}
\end{equation}
holds.
\end{prop}

{\em Proof}: If $n\geq 2$, by (\ref{m-delta}) and (b2) of Definition
\ref{weak-H-mod} we have:
$$u_{n}\wedge u_{n}=\mu_{A}\co (\varphi_{A}\ot \varphi_{A})\co (H\ot
c_{H,A}\ot A)\co ((\delta_{H}\co m_{H}^{n})\ot \eta_{A}\ot
\eta_{A})=\varphi_{A}\co (m_{H}^{n} \ot (\mu_{A}\co (\eta_{A}\ot
\eta_{A})))=u_{n},
$$
and if $n=1$ the equality follows from
$$u_{1}\wedge u_{1}= \mu_{A}\co (\varphi_{A}\ot \varphi_{A})\co (H\ot
c_{H,A}\ot A)\co (\delta_{H}\ot \eta_{A}\ot \eta_{A})=\varphi_{A}\co
(H\ot (\mu_{A}\co (\eta_{A}\ot \eta_{A})))=u_{1}.$$

\begin{defin}
 \label{reg-n}
 {\rm  Let $H$ be a cocommutative weak Hopf algebra and
$(A,\varphi_{A})$ be a weak left $H$-module algebra. For $n\geq 1$,
with
$$Reg_{\varphi_{A}}(H^{n},A)$$ we will denote the set of morphisms
$\sigma:H^n\rightarrow A$ such that there exists a morphism
$\sigma^{-1}:H^n\rightarrow A$ (the convolution inverse of $\sigma$)
satisfying the following equalities:
\begin{itemize}
\item[(c1)] $\sigma\wedge \sigma^{-1}= \sigma^{-1}\wedge
\sigma=u_{n}.$
\item[(c2)] $\sigma\wedge \sigma^{-1}\wedge \sigma=\sigma.$
\item[(c3)] $\sigma^{-1}\wedge \sigma\wedge \sigma^{-1}=\sigma^{-1}.$
\end{itemize}

By $Reg_{\varphi_{A}}(H_{L},A)$ we denote the set of morphisms
$g:H_{L}\rightarrow A$ such that there exists a morphism
$g^{-1}:H_{L}\rightarrow A$ (the convolution inverse of $g$)
satisfying
$$g\wedge g^{-1}=g^{-1}\wedge g=u_{0},\;\; g\wedge g^{-1}
\wedge g=g,\;\; g^{-1}\wedge g\wedge g^{-1}=g^{-1}$$  where
$u_{0}=u_{1}\co i_{L}$. Then, by (b7) of the definition of weak
$H$-module algebra, we have  $u_{1}=u_{0}\co p_{L}$.

Note that the equality
\begin{equation}
\label{u1-reg} \mu_{A}\co (u_{1}\ot \sigma)\co (\delta_{H}\ot
H^{n-1})=\sigma
\end{equation}
holds for all $\sigma\in Reg_{\varphi_{A}}(H^{n}, A)$. Indeed,

\begin{itemize}

\item[ ]$\hspace{0.38cm}\mu_{A}\co (u_{1}\ot \sigma)\co (\delta_{H}\ot
H^{n-1}) $

\item [ ]$=\mu_{A}\co (u_{1}\ot (u_{n}\wedge\sigma))\co (\delta_{H}\ot
H^{n-1}) $

\item [ ]$= \mu_{A}\co ((\mu_{A}\co (u_{1}\ot u_{n}))\ot \sigma)\co (H\ot \delta_{H^{n}})\co
(\delta_{H}\ot H^{n-1})$

\item [ ]$=\mu_{A}\co ((\mu_{A}\co (u_{1}\ot u_{n})\co (\delta_{H}\ot H^{n-1}))\ot \sigma)\co
\delta_{H^n}$

\item [ ]$=u_{n}\wedge\sigma$

\item [ ]$=\sigma$

\end{itemize}

because by (b4) and (b2) of Definition \ref{weak-H-mod} we have
\begin{equation}
\label{u1-reg-1} \mu_{A}\co (u_{1}\ot u_{n})\co (\delta_{H}\ot
H^{n-1})=u_{n}.
\end{equation}
}
\end{defin}

\begin{prop}
\label{new-properties-sigma} Let $H$ be a cocommutative weak Hopf
algebra and let $(A,\varphi_{A})$ be a weak left $H$-module algebra.
Then, for all $\sigma\in Reg_{\varphi_{A}}(H^{n+1},A)$ the following
equalities hold:

\begin{itemize}
\item[(i)] $\sigma\co (H^i\ot ((\Pi_{H}^{L}\ot H)\co \delta_{H})\ot H^{n-i-1})=
\sigma\co (H^{i}\ot \eta_{H}\ot H^{n-i})$ for all $i\in \{0,\dots,
n-1\}.$
\item[(ii)] $\sigma\co (H^{n-1}\ot ((H\ot \Pi_{H}^{R})\co \delta_{H}))=\sigma\co (H^{n}\ot \eta_{H})$
\end{itemize}
\end{prop}

{\em Proof}: First note that if $\sigma\in
Reg_{\varphi_{A}}(H^{n+1},A)$, by (iv) of Proposition
\ref{cocommutative-properties-2} and the equality $\Pi_{H}^{L}\wedge
id_{H}=id_{H}$, we obtain that $\sigma\co (H^i\ot ((\Pi_{H}^{L}\ot
H)\co \delta_{H})\ot H^{n-i-1})\in Reg_{\varphi_{A}}(H^{n},A)$ with
inverse $\sigma^{-1}\co (H^i\ot ((\Pi_{H}^{L}\ot H)\co
\delta_{H})\ot H^{n-i-1})$.

Moreover, by the naturally of $c$ and the equality
(\ref{delta-pi2}), we obtain (i) because:

\begin{itemize}

\item[ ]$\hspace{0.38cm}\sigma\co (H^{i}\ot \eta_{H}\ot H^{n-i})$

\item [ ]$=(u_{n+1}\wedge \sigma)\co (H^{i}\ot \eta_{H}\ot H^{n-i}) $

\item [ ]$=\mu_{A}\co (u_{n}\otimes \sigma)\co (H^{i}\ot \mu_{H}\ot c_{H^{n-i-1},H^{i}}\ot H\ot
 H\ot H^{n-i-1})$
\item[ ]$\hspace{0.38cm} \co (H^i\ot H\ot c_{H^{i},H}\ot c_{H,H^{n-i-1}}\ot H\ot H^{n-i-1})\co
(H^i\ot  c_{H^{i},H}\ot c_{H,H}\ot c_{H,H^{n-i-1}}\ot H^{n-i-1}) $
\item[ ]$\hspace{0.38cm} \co (\delta_{H^{i}}\ot (\delta_{H}\co
\eta_{H})\ot \delta_{H}\ot \delta_{H^{n-i-1}})$

\item [ ]$=\mu_{A}\co (u_{n}\ot \sigma)\co (H^{i}\ot H\ot c_{H^{i}, H^{n-i-1}}\ot H\ot
 H\ot H^{n-i-1})$
\item[ ]$\hspace{0.38cm} \co (H^i \ot c_{H^{i},H}\ot
((c_{H,H^{n-i-1}} \ot H)\co (H\ot c_{H,H^{n-i-1}})\co
(((\Pi_{H}^{L}\ot H)\co \delta_{H})\ot H^{n-i-1})))$
\item[ ]$\hspace{0.38cm} \ot H^{n-i-1})\co (\delta_{H^{i}}\ot \delta_{H}\ot \delta_{H^{n-i-1}})$

\item [ ]$=u_{n}\wedge (\sigma\co (H^i\ot ((\Pi_{H}^{L}\ot H)\co \delta_{H})\ot H^{n-i-1}))$

\item [ ]$=\sigma\co (H^i\ot ((\Pi_{H}^{L}\ot H)\co \delta_{H})\ot H^{n-i-1}), $

\end{itemize}

The proof for (ii) is similar using (\ref{delta-pi4}) and we leave
the details to the reader.

\begin{apart}
{\rm  Let $H$ be a cocommutative weak Hopf algebra and
$(A,\varphi_{A})$ be a weak left $H$-module algebra. Then, $u_{0}\in
Reg_{\varphi_{A}}(H_{L},A)$, $u_{n}\in Reg_{\varphi_{A}}(H^{n},A)$
and $Reg_{\varphi_{A}}(H_{L},A)$, $Reg_{\varphi_{A}}(H^{n},A)$ are
groups with neutral elements $u_{0}$ and $u_{n}$ respectively. Also,
if $A$ is commutative, we have that $Reg_{\varphi_{A}}(H_{L},A)$,
$Reg_{\varphi_{A}}(H^{n},A)$ are abelian groups.

If $(A,\varphi_{A})$ is a left $H$-module algebra, the groups
$Reg_{\varphi_{A}}(H_{L},A)$, $Reg_{\varphi_{A}}(H^{n}, A)$, $n \geq
1$ are the objects of a cosimplicial complex of groups with coface
operators defined by
$$\partial_{0,i}:Reg_{\varphi_{A}}(H_{L},A)\rightarrow Reg_{\varphi_{A}}(H,A), \;\; i\in\{0,1\}$$
$$\partial_{0,0}(g)=\varphi_{A}\co (H\ot (g\co p_{L}\co \Pi_{H}^{R}))\co \delta_{H},\;\;
\partial_{0,1}(g)=g\co p_{L}$$
$$\partial_{1,i}:Reg_{\varphi_{A}}(H,A)\rightarrow Reg_{\varphi_{A}}(H^2,A), \;\; i\in\{0,1,2\}$$
$$\partial_{1,0}(h)=\varphi_{A}\co (H\ot h),\;\;
\partial_{1,1}(h)=h\co \mu_{H},\;\; \partial_{1,2}(h)=h\co \mu_{H}\co (H\ot \Pi_{H}^{L});$$
$$\partial_{k-1,i}:Reg_{\varphi_{A}}(H^{k-1},A)\rightarrow Reg_{\varphi_{A}}(H^{k},A), \;\;k> 2,\;\; i\in\{0,1,\cdots,k\}$$
$$ \partial_{k-1,i}(\sigma)=\left\{ \begin{array}{l}
 \partial_{k-1,0}(\sigma)=\varphi_{A}\co (H\ot \sigma),\\
 \\
 \partial_{k-1,i}(\sigma)=\sigma\co (H^{i-1}\ot\mu_{H}\ot H^{k-i-1}),
 \;\;i\in\{1,\cdots,k-1\} \\

 \\
\partial_{k-1,k}(\sigma)=\sigma\co (H^{k-2}\ot (\mu_{H}\co (H\ot \Pi_{H}^{L}))) ,
  \end{array}\right.
$$

and  codegeneracy operators are defined by
$$s_{1,0}:Reg_{\varphi_{A}}(H,A)\rightarrow Reg_
{\varphi_{A}}(H_{L},A),$$
$$s_{1,0}(h)=h\co i_{L}, $$
$$s_{2,i}:Reg_{\varphi_{A}}(H^{2},A)\rightarrow Reg_
{\varphi_{A}}(H,A),\;\; i\in\{0,1\}$$
$$s_{2,0}(\sigma)=\sigma\co (\eta_{H}\ot H),\;\;\;\;s_{2,1}(\sigma)=\sigma\co (H\ot \eta_{H}),
$$
 and
$$s_{k+1,i}:Reg_{\varphi_{A}}(H^{k+1},A)\rightarrow Reg_
{\varphi_{A}}(H^{k},A), \;\;k\geq 2,\;\; i\in\{0,1,\cdots,k\}$$
\vspace{0.05cm}
$$ s_{k+1,i}(\sigma)=\left\{ \begin{array}{l}
 s_{k+1,0}(\sigma)=\sigma\co( \eta_{H}\ot H^{k}) ,\\
 \\
s_{k+1,i}(\sigma)=\sigma\co (H^{i}\ot \eta_{H}\ot H^{k-i}), \;\;\; i\in\{1,\cdots,k-1\} \\
 \\
s_{k+1,k}(\sigma)=\sigma\co (H^{k}\ot \eta_{H}).
  \end{array}\right.
$$

The morphism $\partial_{0,0}$, is a well defined group morphism because:

\begin{itemize}

\item[ ]$\hspace{0.38cm}\partial_{0,0}(g)\wedge \partial_{0,0}(f)$

\item[ ]$=\mu_{A}\co ((\varphi_{A}\co (H\ot (g\co p_{L}\co \Pi_{H}^{R})))\ot (\varphi_{A}\co (H\ot (f\co p_{L}\co \Pi_{H}^{R})))) \co \delta_{H^{2}}\co \delta_{H}$

\item [ ]$=\mu_{A}\co (\varphi_{A}\ot \varphi_{A})\co (H\ot c_{H,A}\ot A)\co (\delta_{H}\ot (((g\co p_{L}\co \Pi_{H}^{R})\ot (f\co p_{L}\co \Pi_{H}^{R}))\co \delta_{H}))\co \delta_{H}$

\item [ ]$= \varphi_{A}\co (H\ot ((g\co p_{L}\co \Pi_{H}^{R})\wedge (f\co p_{L}\co \Pi_{H}^{R})))\co \delta_{H}$

\item [ ]$=\varphi_{A}\co (H\ot (((g\co p_{L})\wedge (f\co p_{L}))\co \Pi_{H}^{R}) )\co \delta_{H}.$

\item [ ]$=\partial_{0,0}(g\wedge f).$

\end{itemize}
where the first equality follows by (i) of Proposition \ref{cocommutative-properties-2},
the second one by the naturality of $c$, the third one by (b2) of Definition
\ref{weak-H-mod}, the fourth one by (i) of Proposition \ref{cocommutative-properties}
and in the last one was used that $p_{L}$ is a coalgebra morphism (see Remark
\ref{coalgebra-structure}).

Using hat $p_{L}$ is a coalgebra morphism, we obtain that
$\partial_{0,1}$ is a group morphism. Moreover, by (b2) of
Definition \ref{weak-H-mod}, (a1) of Definition \ref{wha},
Proposition \ref{delta-mu-it} and (i) of Proposition
\ref{cocommutative-properties}, we have that $\partial_{k-1,i}$ are
well defined group morphisms for $k\geq 1$.

On the other hand, by (i) of Proposition
\ref{cocommutative-properties} we have that $s_{1,0}$ is a group
morphism and by Propositions \ref{cocommutative-properties} and
\ref{new-properties-sigma} we obtain that $s_{k+1,i}$ are well
defined group morphisms for $k\geq 0$.

We have the cosimplicial identities from the following: For $j=1$,
by (iv) of Proposition \ref{cocommutative-properties} and the
condition of left $H$-module algebra for $A$,  we have
$$\partial_{1,1}(\partial_{0,0}(g))=\varphi_{A}\co (\mu_{H}\ot (g\co
p_{L}\co \Pi_{H}^{R}))\co (H\ot
\delta_{H})=\partial_{1,0}(\partial_{0,0}(g)).$$

Moreover, if $H$ is cocommutative,
$\Pi_{H}^{L}=\overline{\Pi}_{H}^{L}$ and as a consequence
$\Pi_{H}^{L}\co \Pi_{H}^{R}=\Pi_{H}^{L}$. Then, by (i) and (iv) of
Proposition \ref{cocommutative-properties} and the properties of
left $H$-module algebra we get

\begin{itemize}

\item[ ]$\hspace{0.38cm}\partial_{1,2}(\partial_{0,0}(g)) $

\item [ ]$=\varphi_{A}\co (\mu_{H}\ot (g\co
p_{L}\co \Pi_{H}^{R}))\co (H\ot (\delta_{H}\co \Pi_{H}^{L})) $

\item [ ]$=\varphi_{A}\co (\mu_{H}\ot (g\co
p_{L}))\co (H\ot (\delta_{H}\co \Pi_{H}^{L})) $

\item [ ]$=\varphi_{A}\co (H \ot (\varphi_{A}\co (H\ot (g\co
p_{L}))\co \delta_{H}))\co  (H\ot \Pi_{H}^{L}) $

\item [ ]$= \varphi_{A}\co (H \ot (\varphi_{A}\co (\Pi_{H}^{L}\ot (g\co
p_{L}))\co \delta_{H}))$

\item [ ]$= \varphi_{A}\co (H \ot (\mu_{A}\co (u_{1}\ot (g\co
p_{L}))\co \delta_{H}))$

\item [ ]$= \varphi_{A}\co (H \ot (\mu_{A}\co ((u_{0}\co p_{L})\ot (g\co
p_{L}))\co \delta_{H}))$

\item [ ]$= \varphi_{A}\co (H \ot ( (u_{0}\wedge g
)\co p_{L}))$

\item [ ]$= \varphi_{A}\co (H \ot (g
\co p_{L}))$

\item [ ]$= \partial_{1,0}(\partial_{0,1}(g)).$

\end{itemize}

Also, by (\ref{deltamu}) we obtain that
$\partial_{1,2}(\partial_{0,1}(g))=\partial_{1,1}(\partial_{1,0}(g))$.
In a similar way, by the associativity of $\mu_{H}$,
$$\partial_{k,j}\co \partial_{k-1, i}=\partial_{k,i}\co \partial_{k-1, j-1}, \;\;\; j>i $$
for $k>1$.

On the other hand, trivially
$$s_{k-1,j}\co s_{k,i}=s_{k-1,i}\co s_{k,j+1}, \;\;\; j\geq i. $$

Moreover, $s_{1,0}(\partial_{0,0}(g))=\varphi_{A}\co
((\Pi_{H}^{L}\co i_{L})\ot (g\co p_{L}\co \Pi_{H}^{R}\co i_{L}))\co
\delta_{H_{L}}=u_{0}\wedge g=g, $ and
$s_{1,0}(\partial_{0,1}(g))=g$. Also,
$s_{2,0}(\partial_{1,0}(h))=h=s_{2,0}(\partial_{1,1}(h))$,
$s_{2,0}(\partial_{1,2}(h))=h\co
\Pi_{H}^{L}=\partial_{0,1}(s_{1,0}(h)), $
$$s_{2,1}(\partial_{1,0}(h))=\varphi_{A}\co (H\ot (h\co
\Pi_{H}^{R}))\co \delta_{H}=\varphi_{A}\co (H\ot (h\co
\Pi_{H}^{L}\co \Pi_{H}^{R}))\co
\delta_{H}=\partial_{0,0}(s_{1,0}(h))$$ and
$s_{2,1}(\partial_{1,1}(h))=h=s_{2,1}(\partial_{1,2}(h))$ because
$\Pi_{H}^{L}\co \eta_{H}=\eta_{H}$.

Finally, for $k>2$, the identities
$$ s_{k+1,j}\co\partial_{k,i}=\left\{ \begin{array}{l}
\partial_{k-1,i} \co s_{k,j-1},\;\;\; i<j\\
 \\
id_{Reg_{\varphi_{A}}(H^{k},A)}, \;\;\; i=j,\; i=j+1 \\
\\
\partial_{k-1,i-1} \co s_{k,j},\;\;\; i>j+1
  \end{array}\right.
$$
follow as in the Hopf algebra setting.

Let
$$D^{k}_{\varphi_{A}}=\partial_{k,0}\wedge \partial_{k,1}^{-1}\wedge \cdots
\wedge \partial_{k,k+1}^{{(-1)}^{k+1}}$$ be the  coboundary
morphisms of the cochain complex
$$Reg_{\varphi_{A}}(H_{L},A)\stackrel{D^{0}_{\varphi_{A}}}\longrightarrow
Reg_{\varphi_{A}}(H,A)\stackrel{D^{1}_{\varphi_{A}}}\longrightarrow
Reg_{\varphi_{A}}(H^2,A)
\stackrel{D^{2}_{\varphi_{A}}}\longrightarrow \cdots $$
$$
\cdots \stackrel{D^{k-1}_{\varphi_{A}}}\longrightarrow
Reg_{\varphi_{A}}(H^{k},A)\stackrel{D^{k}_{\varphi_{A}}}\longrightarrow
Reg_{\varphi_{A}}(H^{k+1},A)
\stackrel{D^{k+1}_{\varphi_{A}}}\longrightarrow \cdots$$

associated to the cosimplicial complex
$Reg_{\varphi_{A}}(H^{\bullet}, A)$.

Then, when $(A, \varphi_{A})$ is  a commutative left $H$-module
algebra,
$(Reg_{\varphi_{A}}(H^{\bullet},A),D^{\bullet}_{\varphi_{A}})$ gives
the Sweedler cohomology of $H$ in $(A,\varphi_{A})$. Therefore,  the
kth group will be defined by
$$\displaystyle
\frac{Ker(D^{k}_{\varphi_{A}})}{Im(D^{k-1}_{\varphi_{A}})}$$ for
$k\geq 1$ and $Ker(D^{0}_{\varphi_{A}})$ for $k=0$. We will denote
 it by $H^{k}_{\varphi_{A}}(H,A)$.

 The normalized cochain subcomplex   of
$(Reg_{\varphi_{A}}(H^{\bullet},A),D^{\bullet}_{\varphi_{A}})$ is
defined by
$$Reg_{\varphi_{A}}^{+}(H^{k+1}, A)=\bigcap_{i=0}^{k}Ker(s_{k+1,i}),$$
$$Reg_{\varphi_{A}}^{+}(H_{L}, A)=\{g\in Reg_{\varphi_{A}}(H_{L}, A)\;;\; g\co p_{L}\co
\eta_{H}=\eta_{A}\}$$

and $D^{k +}_{\varphi_{A}}$ the restriction of $D^{k}_{\varphi_{A}}$
to $Reg_{\varphi_{A}}^{+}(H^{\bullet}, A)$.

We have that $(Reg_{\varphi_{A}}^{+}(H^{\bullet}, A), D^{\bullet
+}_{\varphi_{A}})$,  is a subcomplex of
$(Reg_{\varphi_{A}}(H^{\bullet},A),D^{\bullet}_{\varphi_{A}})$ and
the injection map induces an isomorphism of cohomology (see
\cite{Maclane} for the dual result). Then,
$$H^{2}_{\varphi_{A}}(H,A)\simeq H^{2+}_{\varphi_{A}}(H,A)= \displaystyle
\frac{Ker(D^{2+}_{\varphi_{A}})}{Im(D^{1+}_{\varphi_{A}})}.$$

Note that
$$Reg_{\varphi_{A}}^{+}(H, A)=Ker(s_{1,0})=\{h \in Reg_{\varphi_{A}}(H, A)\; ;
\; h\co i_{L}=u_{0}\},$$

and
$$Reg_{\varphi_{A}}^{+}(H^{2}, A)=Ker(s_{2,0})\cap Ker(s_{2,1})$$
$$=\{\sigma \in Reg_{\varphi_{A}}(H^{2}, A)\;; \; \sigma\co
(\eta_{H}\ot H)=\sigma\co (H\ot \eta_{H})=u_{1}\}.$$

The following proposition give a different characterization
of the morphisms in $Reg^{+}_{\varphi_{A}}(H,A)$.}
\end{apart}

\begin{prop}
\label{f-eta} Let $H$ be a weak Hopf algebra and $(A,\varphi_{A})$
be a weak left $H$-module algebra. Let $h:H\rightarrow A$ be a
morphism satisfying
$$ h\wedge h^{-1}=
h^{-1}\wedge h=u_{1}, \;\; h\wedge h^{-1}\wedge h=h,
\;\;h^{-1}\wedge h\wedge h^{-1}=h^{-1}.$$ The following equalities
are equivalent
\begin{itemize}

\item[(i)] $h\co \eta_{H}=\eta_{A}.$
\item[(ii)] $h\co \Pi_{H}^{L}=u_{1}.$
\item[(iii)] $h\co \overline{\Pi}_{H}^{L}=u_{1}.$
\end{itemize}

\end{prop}

{\em Proof}: The assertion (ii) $\Rightarrow$ (i) follows by
$$h\co \eta_{H}=h\co \Pi_{H}^{L}\co \eta_{H}=u_{1}\co
\eta_{H}=\eta_{A}.$$

Now we get (i) $\Rightarrow$ (ii):

\begin{itemize}

\item[ ]$\hspace{0.38cm}h\co \Pi_{H}^{L} $

\item [ ]$=(u_{1}\wedge  h)\co \Pi_{H}^{L} $

\item [ ]$=((\varepsilon_{H}\co \mu_{H})\ot (\mu_{A}\co (u_{1}\ot h))\co
(H\ot c_{H,H}\ot H)\co (\delta_{H}\ot c_{H,H})\co ((\delta_{H}\co
\eta_{H})\ot H) $

\item [ ]$= \mu_{A}\co (u_{2}\ot h)\co (H\ot c_{H,H})\co
((\delta_{H}\co \eta_{H})\ot\Pi_{H}^{L}) $

\item [ ]$= \mu_{A}\co (\varphi_{A}\ot h)\co (H\ot c_{H,A})\co
((\delta_{H}\co \eta_{H})\ot (u_{1}\co \Pi_{H}^{L}))$

\item [ ]$=\mu_{A}\co (\varphi_{A}\ot h)\co (H\ot c_{H,A})\co
((\delta_{H}\co \eta_{H})\ot u_{1}) $

\item [ ]$=\mu_{A}\co (\varphi_{A}\ot h)\co (\overline{\Pi}_{H}^{L}\ot c_{H,A})\co
((\delta_{H}\co \eta_{H})\ot u_{1}) $

\item [ ]$=\mu_{A}\co ((\mu_{A}\co c_{A,A}\co (u_{1}\ot A))\ot h)\co  (H\ot c_{H,A})\co
((\delta_{H}\co \eta_{H})\ot u_{1})$

\item [ ]$= \mu_{A}\co c_{A,A}\co (((u_{1}\wedge h)\co\eta_{H})\ot u_{1})$

\item [ ]$=\mu_{A}\co c_{A,A}\co (( h\co\eta_{H})\ot u_{1}) $

\item [ ]$=u_{1}. $

\end{itemize}

The first equality follows by the properties of $h$, the second one
by the naturality of $c$ and the coassociativity of $\delta_{H}$,
the third one by (\ref{delta-pi1}), the fourth one by (b3) of
Definition \ref{weak-H-mod}, the fifth one by (b6) of Definition
\ref{weak-H-mod}, the sixth one by (\ref{delta-pi31}), the seventh
one by (b5) of Definition \ref{weak-H-mod}, the eight one by the
naturality of $c$ and the associativity of $\mu_{A}$, the ninth one
the by the properties of $h$ and the last one by (ii).

The assertion (iii) $\Rightarrow$ (i) follows because
$$h\co \eta_{H}=h\co \overline{\Pi}_{H}^{L}\co \eta_{H}=u_{1}\co
\eta_{H}=\eta_{A}.$$

The proof for (i) $\Rightarrow$ (iii) is the following:

\begin{itemize}

\item[ ]$\hspace{0.38cm}h\co \overline{\Pi}_{H}^{L} $

\item [ ]$=(u_{1}\wedge  h)\co  \overline{\Pi}_{H}^{L} $

\item [ ]$= \mu_{A}\co (h\ot ((u_{1}\ot (\varepsilon_{H}\co \mu_{H}))\co
(\delta_{H}\ot H)))\co ((\delta_{H}\co \eta_{H})\ot H)$

\item [ ]$=\mu_{A}\co (h\ot (u_{1}\co\mu_{H}\co ( \overline{\Pi}_{H}^{L}\ot H)))\co
((\delta_{H}\co \eta_{H})\ot H) $

\item [ ]$=\mu_{A}\co (h\ot\varphi_{A})\co
((\delta_{H}\co \eta_{H})\ot u_{1}) $

\item [ ]$=\mu_{A}\co (h\ot(\varphi_{A}\co (\Pi_{H}^{L}\ot A)))\co
((\delta_{H}\co \eta_{H})\ot u_{1}) $

\item [ ]$= \mu_{A}\co (h\ot(\mu_{A}\co (u_{1}\ot A)))\co
((\delta_{H}\co \eta_{H})\ot u_{1})$

\item [ ]$=\mu_{A}\co (((h\wedge u_{1})\co \eta_{H})\ot u_{1}) $

\item [ ]$=\mu_{A}\co ((h\co \eta_{H})\ot u_{1}) $

\item [ ]$=u_{1}. $

\end{itemize}

The first equality follows by the properties of $h$, the second one
by  the coassociativity of $\delta_{H}$, the third one by
(\ref{delta-pi21}), the fourth one by (b3) of Definition
\ref{weak-H-mod}, the fifth one by (b7) of Definition
\ref{weak-H-mod} and (\ref{delta-pi2}), the sixth one by (b4) of
Definition \ref{weak-H-mod}, the seventh one by the associativity of
$\mu_{H}$, the eight one by the properties of $h$ and the last one
by (iii).

\begin{rem}
{\rm
 Note that as a consequence of  Proposition \ref{f-eta}:
$$Reg_{\varphi_{A}}^{+}(H, A)=\{h \in Reg_{\varphi_{A}}(H,
A)\; ; \; h\co \eta_{H}=\eta_{A}\},$$ and by Proposition
\ref{new-properties-sigma}  we have
$$Reg_{\varphi_{A}}^{+}(H^{2}, A)=\{\sigma \in Reg_{\varphi_{A}}(H^{2}, A)\; ; \; \sigma\co
(\Pi_{H}^{L}\ot H)\co \delta_{H}=\sigma\co (H\ot \Pi_{H}^{R})\co
\delta_{H}=u_{1}\}.$$

}
\end{rem}

\section{Weak crossed products for weak Hopf algebras}

In the first paragraphs of this section we resume some basic facts
about the general theory of weak crossed products in  $\mathcal C$
introduced in \cite{mra-preunit} particularized for a weak Hopf
algebra $H$.

Let $A$ be an algebra and let $H$ be a weak Hopf algebra in
${\mathcal C}$. Suppose that there exists a morphism
$$\psi_{H}^{A}:H\ot A\rightarrow A\ot H$$  such that the following
equality holds
\begin{equation}\label{wmeas-wcp}
(\mu_A\ot H)\co (A\ot \psi_{H}^{A})\co (\psi_{H}^{A}\ot A) =
\psi_{H}^{A}\co (H\ot \mu_A).
\end{equation}

 As a consequence of (\ref{wmeas-wcp}), the morphism $\nabla_{A\ot H}:A\ot H\rightarrow
A\ot H$ defined by
\begin{equation}\label{idem-wcp}
\nabla_{A\ot H} = (\mu_A\ot H)\co(A\ot \psi_{H}^{A})\co (A\ot H\ot
\eta_A)
\end{equation}
is an idempotent. Moreover, it satisfies that
$$\nabla_{A\ot H}\co (\mu_A\ot H) = (\mu_A\ot H)\co
(A\ot \nabla_{A\ot H}),$$ that is, $\nabla_{A\ot H}$ is a left
$A$-module morphism (see Lemma 3.1 of \cite{mra-preunit}) for the
regular action  $\varphi_{A\ot H}=\mu_{A}\ot H$. With $A\times H$,
$i_{A\ot H}:A\times H\rightarrow A\ot H$ and $p_{A\ot H}:A\ot
H\rightarrow A\times H$ we denote the object, the injection and the
projection associated to the factorization of $\nabla_{A\ot H}$.
Finally, if $\psi_{H}^{A}$ satisfies (\ref{wmeas-wcp}), the
following identities hold
\begin{equation}\label{fi-nab}
(\mu_{A}\ot H)\co (A\ot \psi_{H}^{A})\co (\nabla_{A\ot H}\ot A)=
(\mu_{A}\ot H)\co (A\ot \psi_{H}^{A})=\nabla_{A\ot H}\co(\mu_{A}\ot
H)\co (A\ot \psi_{H}^{A}).
\end{equation}

From now on we consider quadruples ${\Bbb A}_{H}=(A, H,
\psi_{H}^{A}, \sigma_{H}^{A})$ where $A$ is an algebra, $H$ an
object, $\psi_{H}^{A}:H\ot A\rightarrow A\ot H$ a morphism
satisfiying (\ref{wmeas-wcp}) and $\sigma_{H}^{A}:H\ot H\rightarrow
A\ot H$  a morphism in ${\mathcal C}$.

We say that ${\Bbb A}_{H}=(A, H, \psi_{H}^{A}, \sigma_{H}^{A})$
satisfies the twisted condition if
\begin{equation}\label{twis-wcp}
(\mu_A\ot H)\co (A\ot \psi_{H}^{A})\co (\sigma_{H}^{A}\ot A) =
(\mu_A\ot H)\co (A\ot \sigma_{H}^{A})\co (\psi_{H}^{A}\ot H)\co
(H\ot \psi_{H}^{A})
\end{equation}
and   the  cocycle condition holds if
\begin{equation}\label{cocy2-wcp}
(\mu_A\ot H)\co (A\ot \sigma_{H}^{A}) \co (\sigma_{H}^{A}\ot H) =
(\mu_A\ot H)\co (A\ot \sigma_{H}^{A})\co (\psi_{H}^{A}\ot H)\co
(H\ot\sigma_{H}^{A}).
\end{equation}

Note that, if ${\Bbb A}_{H}=(A, H, \psi_{H}^{A}, \sigma_{H}^{A})$
satisfies the twisted condition, in Proposition 3.4 of
\cite{mra-preunit} we prove that the following equalities hold:
\begin{equation}\label{c1}
(\mu_A\otimes H)\circ (A\otimes \sigma_{H}^{A})\circ
(\psi_{H}^{A}\otimes H)\circ (H\otimes \nabla_{A\otimes H}) =
\nabla_{A\otimes H}\circ (\mu_A\otimes H)\circ (A\otimes
\sigma_{H}^{A})\circ (\psi_{H}^{A}\otimes H),
\end{equation}
\begin{equation}\label{aw}
\nabla_{A\otimes H}\circ (\mu_A\otimes H)\circ
(A\otimes\sigma_{H}^{A})\circ (\nabla_{A\otimes H}\otimes H) =
\nabla_{A\otimes H}\circ (\mu_A\otimes H)\circ
(A\otimes\sigma_{H}^{A}).
\end{equation}

Then, if $\nabla_{A\ot H}\co\sigma_{H}^{A}=\sigma_{H}^{A}$ we obtain
\begin{equation}\label{c11}
(\mu_A\otimes H)\circ (A\otimes \sigma_{H}^{A})\circ
(\psi_{H}^{A}\otimes H)\circ (H\otimes \nabla_{A\otimes H}) =
 (\mu_A\otimes H)\circ (A\otimes
\sigma_{H}^{A})\circ (\psi_{H}^{A}\otimes H),
\end{equation}
\begin{equation}\label{aw1}
 (\mu_A\otimes H)\circ
(A\otimes\sigma_{H}^{A})\circ (\nabla_{A\otimes H}\otimes H) =
(\mu_A\otimes H)\circ (A\otimes\sigma_{H}^{A}).
\end{equation}

By virtue of (\ref{twis-wcp}) and (\ref{cocy2-wcp}) we will consider
from now on, and without loss of generality, that
\begin{equation}
\label{idemp-sigma-inv} \nabla_{A\ot H}\co\sigma_{H}^{A} =
\sigma_{H}^{A}
\end{equation}

holds for all quadruples ${\Bbb A}_{H}=(A, H, \psi_{H}^{A},
\sigma_{H}^{A})$ {(see Proposition 3.7 of \cite{mra-preunit})}.

For ${\Bbb A}_{H}=(A, H, \psi_{H}^{A}, \sigma_{H}^{A})$ define the
product
\begin{equation}\label{prod-todo-wcp}
\mu_{A\ot  H} = (\mu_A\ot H)\co (\mu_A\ot \sigma_{H}^{A})\co (A\ot
\psi_{H}^{A}\ot H)
\end{equation}
and let $\mu_{A\times H}$ be the restriction of $\mu_{A\ot  H}$ to
$A\times H$, i.e.
\begin{equation}
\label{prod-wcp} \mu_{A\times H} = p_{A\ot H}\co\mu_{A\ot H}\co
(i_{A\ot H}\ot i_{A\ot H}).
\end{equation}

If the twisted and the cocycle conditions hold, the product
$\mu_{A\ot H}$ is associative and normalized with respect to
$\nabla_{A\ot H}$ (i.e. $\nabla_{A\ot H}\co \mu_{A\ot H}=\mu_{A\ot
H}=\mu_{A\ot H}\co (\nabla_{A\ot H}\ot \nabla_{A\ot H}$)) and, by
the definition of $\mu_{A\ot H}$,
\begin{equation}
\label{otra-prop} \mu_{A\ot  H}\co (\nabla_{A\ot H}\ot A\ot
H)=\mu_{A\ot H}
\end{equation}
holds and therefore
\begin{equation}
\label{vieja-proof} \mu_{A\ot  H}\circ (A\otimes H\otimes
\nabla_{A\otimes H})=\mu_{A\otimes H}.
\end{equation}
Due to the normality condition, $\mu_{A\times H}$ is associative as
well (Proposition 2.5 of \cite{mra-preunit}). Hence we define:
\begin{defin}\label{wcp-def}{\rm
If ${\Bbb A}_{H}=(A,H, \psi_{H}^{A}, \sigma_{H}^{A})$  satisfies
(\ref{twis-wcp}) and (\ref{cocy2-wcp}) we say that $(A\ot H,
\mu_{A\ot H})$ is a weak crossed product.}
\end{defin}

The next natural question that arises is if it is possible to endow
$A\times H$ with a unit, and hence with an algebra structure. As we
recall in \cite{mra-preunit}, in order to do that, we need to use
the notion of preunit to obtain an unit in $A\times H$. In our
setting, if $A$ is an algebra, $H$ an object in ${\mathcal C}$ and
$m_{A\otimes H}$ is an associative product defined in $A\otimes H$ a
preunit $\nu:K\rightarrow A\otimes H$ is a morphism satisfying
\begin{equation}
m_{A\otimes H}\circ (A\otimes H\otimes \nu)=m_{A\otimes H}\circ
(\nu\otimes A\otimes H)=m_{A\otimes H}\circ (A\otimes H\otimes
(m_{A\otimes H}\circ (\nu\otimes \nu))).
\end{equation}

Associated to a preunit we obtain an idempotent morphism
$$\nabla_{A\otimes
H}^{\nu}=m_{A\otimes H}\circ (A\otimes H\otimes \nu):A\otimes
H\rightarrow A\otimes H.$$

 Take $A\times^{\nu} H$ the image of this
idempotent, $p_{A\otimes H}^{\nu}$ the projection and $i_{A\otimes
H}^{\nu}$ the injection. It is possible to endow $A\times^{\nu} H$
with an algebra structure whose product is
$$m_{A\times^{\nu} H} = p_{A\otimes
H}^{\nu}\circ m_{A\otimes H}\circ (i_{A\otimes H}^{\nu}\otimes
i_{A\otimes H}^{\nu})$$ and whose unit is $\eta_{A\times^{\nu}
H}=p_{A\otimes H}^{\nu}\circ \nu$ (see Proposition 2.5 of
\cite{mra-preunit}). If moreover, $\mu_{A\otimes H}$ is left
$A$-linear for the actions $\varphi_{A\otimes H}=\mu_{A}\otimes H$,
$\varphi_{A\otimes H\otimes A\otimes H }=\varphi_{A\otimes H}\otimes
A\otimes H$ and normalized with respect to $\nabla_{A\otimes
H}^{\nu}$,  the morphism
\begin{equation}
\label{beta-nu} \beta_{\nu}:A\rightarrow A\otimes H,\; \beta_{\nu} =
(\mu_A\otimes H)\circ (A\otimes \nu)
\end{equation}

is multiplicative and left $A$-linear for $\varphi_{A}=\mu_{A}$.

Although $\beta_{\nu}$ is not an algebra morphism, because $A\otimes
H$ is not an algebra, we have that $\beta_{\nu}\circ \eta_A = \nu$,
and thus the morphism $\bar{\beta_{\nu}} = p_{A\otimes
H}^{\nu}\circ\beta_{\nu}:A\rightarrow A\times^{\nu} H$ is an algebra
morphism.

In light of the considerations made in the last paragraphs, and
using the twisted and the cocycle conditions, in \cite{mra-preunit}
we characterize weak crossed products with a preunit, and moreover
we obtain an algebra structure on $A\times H$. These assertions are
a consequence of the following theorem proved in \cite{mra-preunit}.

\begin{teo}
\label{thm1-wcp} Let $A$ be an algebra, $H$ a weak Hopf algebra and
$m_{A\otimes H}:A\otimes H\otimes A\otimes H\rightarrow A\otimes H$
a morphism of left $A$-modules  for the actions $\varphi_{A\otimes
H}=\mu_{A}\otimes H$, $\varphi_{A\otimes H\otimes A\otimes H
}=\varphi_{A\otimes H}\otimes  A\otimes H$.

Then the following statements are equivalent:
\begin{itemize}
\item[(i)] The product $m_{A\otimes H}$ is associative with preunit
$\nu$ and normalized with respect to $\nabla_{A\otimes H}^{\nu}.$

\item[(ii)] There exist morphisms $\psi_{H}^{A}:H\otimes A\rightarrow A\otimes
V$, $\sigma_{H}^{A}:H\otimes H\rightarrow A\otimes H$ and
$\nu:k\rightarrow A\otimes H$ such that if $\mu_{A\otimes H}$ is the
product defined in (\ref{prod-todo-wcp}), the pair $(A\otimes H,
\mu_{A\otimes H})$ is a weak crossed product with $m_{A\otimes H} =
\mu_{A\otimes H}$ and satisfying:
\begin{equation}\label{pre1-wcp}
    (\mu_A\otimes H)\circ (A\otimes \sigma_{H}^{A})\circ
    (\psi_{H}^{A}\otimes H)\circ (H\otimes \nu) =
    \nabla_{A\otimes H}\circ
    (\eta_A\otimes H),
    \end{equation}
\begin{equation}\label{pre2-wcp}
    (\mu_A\otimes H)\circ (A\otimes \sigma_{H}^{A})\circ
    (\nu\otimes H) = \nabla_{A\otimes H}\circ (\eta_A\otimes H),
    \end{equation}
\begin{equation}\label{pre3-wcp}
(\mu_A\otimes H)\circ (A\otimes \psi_{H}^{A})\circ (\nu\otimes A) =
\beta_{\nu},
\end{equation}

\end{itemize}
where $\beta_{\nu}$ is the morphism defined in (\ref{beta-nu}). In
this case $\nu$ is a preunit for $\mu_{A\otimes H}$, the idempotent
morphism of the weak crossed product $\nabla_{A\otimes H}$ is the
idempotent $\nabla_{A\otimes H}^{\nu}$, and we say that the pair
$(A\otimes H, \mu_{A\otimes H})$ is a weak crossed product with
preunit $\nu$.
\end{teo}

\begin{rem}
\label{proof-resume} {\rm Note that in the proof of the previous
Theorem for $(i)\;\Rightarrow \;(ii)$ we define $\psi_{H}^{A}$ and
$\sigma_{H}^{A}$ as
\begin{equation}\label{fi-wcp}
\psi_{H}^{A} = m_{A\otimes H}\circ (\eta_A\otimes H\otimes
\beta_{\nu}),
\end{equation}
\begin{equation}\label{sigma-wcp}
\sigma_{H}^{A} = m_{A\otimes H}\circ (\eta_A\otimes H\otimes
\eta_A\otimes H).
\end{equation}

Also, (\ref{pre3-wcp}) implies that $\nabla_{A\otimes H}\co
\nu=\nu$. }

\end{rem}

\begin{cor}\label{corol-wcp}
If $(A\otimes H, \mu_{A\otimes H})$ is a weak crossed product with
preunit $\nu$, then $A\times H$ is an algebra with the product
defined in (\ref{prod-wcp}) and unit $\eta_{A\times H}=p_{A\otimes
H}\circ\nu$.
\end{cor}

\begin{rem}
\label{corol-wcp-1} {\rm As a consequence of the previous corollary
we obtain that if a   weak crossed product $(A\otimes H,
\mu_{A\otimes H})$ admits two preunits $\nu_{1}$, $\nu_{2}$, as in
(ii) of Theorem \ref{thm1-wcp}, we have $p_{A\otimes
H}\circ\nu_{1}=p_{A\otimes H}\circ\nu_{2}$ and then
$$\nu_{1}=\nabla_{A\ot H}\co \nu_{1}=\nabla_{A\ot H}\co \nu_{2}=\nu_{2}.$$
}
\end{rem}

\begin{defin}
\label{psi-sigma} {\rm Let $H$ be a  weak Hopf algebra,
$(A,\varphi_{A})$  a weak left $H$-module algebra and $\sigma:H\ot
H\rightarrow A$ a morphism. We define the morphisms
$$\psi_{H}^{A}:H\ot A\rightarrow A\ot H, \;\;\;\sigma_{H}^{A}:H\ot H\rightarrow A\ot
H, $$ by
\begin{equation}
\label{psiAH} \psi_{H}^{A}=(\varphi_{A}\ot H)\co (H\ot c_{H,A})\co
(\delta_{H}\ot A)
\end{equation}
and
\begin{equation}
\label{sigmaAH} \sigma_{H}^{A}=(\sigma\ot \mu_{H})\co \delta_{H^2}.
\end{equation}

}
\end{defin}

\begin{prop}
\label{psi-prop} Let $H$ be a weak Hopf algebra and
$(A,\varphi_{A})$  a weak left $H$-module algebra. The morphism
$\psi_{H}^{A}$ defined above satisfies (\ref{wmeas-wcp}). As a
consequence  the morphism $\nabla_{A\ot H}$, defined in
(\ref{idem-wcp}), is an idempotent and the following equalities
hold:
\begin{equation}
\label{nabla-nabla} \nabla_{A\ot H}=((\mu_{A}\co (A\ot u_{1}))\ot
H)\co (A\ot \delta_{H}),
\end{equation}
\begin{equation}
\label{nabla-fi} \mu_{A}\co (u_{1}\ot \varphi_{A})\co (\delta_{H}\ot
A)=\varphi_{A},
\end{equation}
\begin{equation}
\label{nabla-fiAH} (\mu_{A}\ot H)\co (u_{1}\ot \psi_{H}^{A})\co
(\delta_{H}\ot A)=\psi_{H}^{A},
\end{equation}
\begin{equation}
\label{eta-psi-varep} (A\ot \varepsilon_{H})\co \psi_{H}^{A}\co
(H\ot \eta_{A})=u_{1},
\end{equation}
\begin{equation}
\label{eta-psi-complex} (\mu_{A}\ot H)\co (u_{1}\ot c_{H,A})\co
(\delta_{H}\ot A)=(\mu_{A}\ot H)\co (A\ot c_{H,A})\co
((\psi_{H}^{A}\co (H\ot \eta_{A}))\ot A),
\end{equation}
\begin{equation}
\label{nabla-varep} (A\ot \varepsilon_{H})\co \nabla_{A\otimes
H}=\mu_{A}\co (A\ot u_{1}).
\end{equation}
\begin{equation}
\label{nabla-delta} (A\ot \delta_{H})\co \nabla_{A\ot
H}=(\nabla_{A\ot H}\ot H)\co (A\ot \delta_{H}).
\end{equation}
\end{prop}
{\em Proof}: For the morphism $\psi_{H}^{A}$ we have
\begin{itemize}

\item[ ]$\hspace{0.38cm} (\mu_A\ot H)\co (A\ot \psi_{H}^{A})\co (\psi_{H}^{A}\ot A)$

\item [ ]$= ((\mu_{A}\co
(\varphi_{A}\ot \varphi_{A})\co (H\ot c_{H,A}\ot A)\co
(\delta_{H}\ot A\ot A))\ot H)\co (H\ot A\ot c_{H,A})\co$
\item[ ]$\hspace{0.38cm}  (H\ot
c_{H,A}\ot A)\co (\delta_{H}\ot A\ot A)$

\item [ ]$= (
(\varphi_{A}\co (H\ot \mu_{A}))\ot H)\co (H\ot A\ot c_{H,A})\co
(H\ot c_{H,A}\ot A)\co (\delta_{H}\ot A\ot A) $

\item [ ]$=\psi_{H}^{A}\co (H\ot \mu_A).$

\end{itemize}

where the first equality follows by the naturality of $c$ and the
coassociativity of $\delta_{H}$, the second one by (b2) of
Definition (\ref{weak-H-mod}) and the third one by  the naturality
of $c$. Thus, $\psi_{H}^{A}$ satisfies (\ref{wmeas-wcp}). As a
consequence, $\nabla_{A\ot H}$ is an idempotent and
(\ref{nabla-nabla}),(\ref{eta-psi-varep}), (\ref{nabla-varep})
follow easily from the definition of $\psi_{H}^{A}$. On the other
hand, (\ref{nabla-fi}) follows by (\ref{nabla-nabla}) and (b2) of
Definition \ref{weak-H-mod} because:

\begin{itemize}

\item[ ]$\hspace{0.38cm}\mu_{A}\co (u_{1}\ot \varphi_{A})\co (\delta_{H}\ot
A) $

\item [ ]$=\mu_{A}\co (A\ot \varphi_{A})\co ((\nabla_{A\ot H}\co
(\eta_{A}\ot H))\ot A) $

\item [ ]$=\mu_{A}\co (\varphi_{A}\ot\varphi_{A})\co
(H\ot c_{H,A}\ot A)\co (\delta_{H}\ot \eta_{A}\ot A)  $

\item [ ]$=\varphi_{A}.$

\end{itemize}

Analogously, by (b2) of Definition \ref{weak-H-mod}, we obtain
(\ref{nabla-fiAH}). Finally, the equality (\ref{nabla-delta})
follows from (\ref{nabla-nabla}) and the coassociativity of
$\delta_{H}$,  and (\ref{eta-psi-complex}) is an easy consequence of
the naturality of $c$.

\begin{prop}
\label{sigma-prop1} Let $H$ be a  weak Hopf algebra,
$(A,\varphi_{A})$  a weak left $H$-module algebra and $\sigma:H\ot
H\rightarrow A$  a morphism. The morphism $\sigma_{H}^{A}$
introduced in Definition \ref{psi-sigma} satisfies the following
identity:
\begin{equation}
\label{delta-sigmaHA} (A\ot \delta_{H})\co
\sigma_{H}^{A}=(\sigma_{H}^{A}\ot \mu_{H})\co \delta_{H^{2}}.
\end{equation}

\end{prop}

{\em Proof}: The proof  is the following:
\begin{itemize}

\item[ ]$\hspace{0.38cm}(A\ot \delta_{H})\co \sigma_{H}^{A} $

\item [ ]$=(A\ot \mu_{H}\ot \mu_{H})\co (\sigma\ot  \delta_{H^{2}})\co \delta_{H^{2}}$

\item [ ]$= (\sigma\ot \mu_{H}\ot \mu_{H})\co (H\ot c_{H,H}\ot c_{H,H}\ot H)\co (H\ot H\ot c_{H,H}\ot H\ot H)\co $
\item[ ]$\hspace{0.38cm}(((\delta_{H}\ot H)\co \delta_{H})\ot ((\delta_{H}\ot H)\co \delta_{H}))$

\item [ ]$= (\sigma_{H}^{A}\ot
\mu_{H})\co \delta_{H^{2}}.$

\end{itemize}

The first equality follows by (a1) of Definition \ref{wha}, the
second one by the naturality of $c$ and the last one by the
coassociativity of $\delta_{H}$ and the  naturality of $c$.

\begin{prop}
\label{sigma-prop} Let $H$ be a cocommutative weak Hopf algebra,
$(A,\varphi_{A})$  a weak left $H$-module algebra and $\sigma\in
Reg_{\varphi_{A}}(H^2,A)$. The morphism $\sigma_{H}^{A}$ introduced
in Definition \ref{psi-sigma} satisfies the following identities:

\begin{itemize}

\item[(i)] $\nabla_{A\ot H}\co\sigma_{H}^{A}=\sigma_{H}^{A}.$

\item[(ii)] $(A\ot \varepsilon_{H})\co \sigma_{H}^{A}=\sigma.$

\end{itemize}
\end{prop}

{\em Proof:} By  Proposition \ref{sigma-prop1} and the
properties of $\sigma$ we have that

\begin{itemize}

\item[ ]$\hspace{0.38cm}\nabla_{A\ot H}\co\sigma_{H}^{A} $

\item [ ]$=((\mu_{A}\co (A\ot u_{1}))\ot H)\co (A\ot \delta_{H})\co \sigma_{H}^{A}  $

\item [ ]$=((\mu_{A}\co (A\ot u_{1}))\ot H)\co (\sigma_{H}^{A}\ot
\mu_{H})\co \delta_{H^{2}} $

\item [ ]$= ((\sigma\wedge \sigma^{-1}\wedge \sigma)\ot \mu_{H})\co \delta_{H^{2}} $

\item[ ]$=\sigma_{H}^{A}, $

\end{itemize}

and therefore, (i) holds. Finally, the proof for (ii) follows by
(\ref{nabla-varep}) and (ii) because:
$$(A\ot \varepsilon_{H})\co \sigma_{H}^{A}=(A\ot \varepsilon_{H})\co
 \nabla_{A\otimes H}\co \sigma_{H}^{A}= \mu_{A}\co (A\ot u_{1})\co \sigma_{H}^{A}=
 \sigma\wedge u_{2}=\sigma.$$

\begin{rem}
\label{quadruple-1} {\rm Let $H$ be a cocommutative weak Hopf
algebra, $(A,\varphi_{A})$  a weak left $H$-module algebra and
$\sigma\in Reg_{\varphi_{A}}(H^2,A)$. Note that, by Propositions
\ref{psi-prop}, \ref{sigma-prop1} and \ref{sigma-prop}, we have a
quadruple ${\Bbb A}_{H}=(A,H, \psi_{H}^{A}, \sigma_{H}^{A})$  such
that $\psi_{H}^{A}$ satisfies (\ref{wmeas-wcp}) and $\nabla_{A\ot
H}\co \sigma_{H}^{A}=\sigma_{H}^{A}.$

}
\end{rem}

\begin{defin}
{\rm Let $H$ be a cocommutative weak Hopf algebra, $(A,\varphi_{A})$
a weak left $H$-module algebra and $\sigma\in
Reg_{\varphi_{A}}(H^2,A)$.  We say that $\sigma$ satisfies the
twisted condition if
\begin{equation}
\label{twisted-sigma} \mu_{A}\co ((\varphi_{A}\co (H\ot
\varphi_{A}))\ot A)\co (H\ot H\ot c_{A,A})\co (((H\ot H\ot
\sigma)\co \delta_{H^{2}}) \ot A)=\mu_{A}\co (A\ot \varphi_{A})\co
(\sigma_{H}^{A}\ot A).
\end{equation}

If
\begin{equation}
\label{2-cocycle-sigma}
\partial_{2,3}(\sigma)\wedge \partial_{2,1}(\sigma)=\partial_{2,0}(\sigma)\wedge
\partial_{2,2}(\sigma)
\end{equation}
 holds, we will say that $\sigma$ satisfies the
2-cocycle condition.

}
\end{defin}

\begin{apart}
{\rm
Let $H$ be a weak Hopf algebra. The morphisms
\begin{equation}
\label{omega-L} \Omega_{H\ot H}^{L}= ((\varepsilon_{H}\co
\mu_{H})\ot H\ot H)\co \delta_{H\ot H}:H\ot H\rightarrow H\ot H
\end{equation}
\begin{equation}
\label{omega-R} \Omega_{H\ot H}^{R}= (H\ot H\ot (\varepsilon_{H}\co
\mu_{H}))\co \delta_{H\ot H}:H\ot H\rightarrow H\ot H
\end{equation}
are idempotent. Indeed: By (\ref{delta-pi1}) we have
\begin{equation}
\label{omega-L-1} \Omega_{H\ot H}^{L}=((\mu_{H}\co (H\ot
\Pi_{H}^{L}))\ot H)\co(H\ot \delta_{H}).
\end{equation}

Then, by (\ref{omega-L-1}), the coassociativity of $\delta_{H}$ and
(\ref{id}) we have
$$\Omega_{H\ot H}^{L}\co \Omega_{H\ot H}^{L}=((\mu_{H}\co (H \ot (\Pi_{H}^{L}\wedge
\Pi_{H}^{L}))\ot H)\co (H\ot \delta_{H})=\Omega_{H\ot H}^{L}.$$

The proof for $\Omega_{H\ot H}^{R}$ is similar, using the identity
\begin{equation}
\label{omega-R-1} \Omega_{H\ot H}^{R}=(H\ot (\mu_{H}\co (\Pi^{R}_{H}
\ot H))\co ( \delta_{H}\ot H),
\end{equation} and we left the details to the reader.

By  (a1) of
Definition \ref{wha} we obtain that
\begin{equation}
\label{omega-mu}
\mu_{H}\co \Omega_{H\ot H}^{L}=\mu_{H}\co \Omega_{H\ot
H}^{R}=\mu_{H}
\end{equation}
and it is easy to show that, if we consider the left-right $H$-module actions and a left-right $H$-comodule
coactions
$$\varphi_{H\ot H}=\mu_{H}\ot H,\;\;\phi_{H\ot H}=H\ot \mu_{H}, \;\;
\varrho_{H\ot H}=\delta_{H}\ot H,\;\; \rho_{H\ot H}=H\ot
\delta_{H}$$ on $H\ot H$, we have that $\Omega_{H\ot H}^{L}$ is a
morphism of left and right $H$-modules and right $H$-comodules and
$\Omega_{H\ot H}^{R}$ is a morphism of  left and right $H$-modules
and left $H$-comodules.  Moreover, if $H$ is cocommutative it is an
easy exercise to prove that $\Omega_{H\ot H}^{L}=\Omega_{H\ot
H}^{R}$ and the following equalities hold:
\begin{equation}
\label{firts-cocom-1} \delta_{H\ot H}\co \Omega_{H\ot H}^{L}=(H\ot
H\ot\Omega_{H\ot H}^{L})\co \delta_{H\ot H}=(\Omega_{H\ot H}^{L}\ot
H\ot H)\co \delta_{H\ot H}.
\end{equation}

As a consequence, we obtain that
\begin{equation}
\label{firts-cocom-2} \delta_{H\ot H}\co \Omega_{H\ot
H}^{L}=(\Omega_{H\ot H}^{L}\ot \Omega_{H\ot H}^{L})\co \delta_{H\ot
H}.
\end{equation}

Therefore, if $H$ is cocommutative, we will denote the morphism $\Omega_{H\ot H}^{L}$ by $\Omega_{H}^{2}$.
}
\end{apart}

\begin{prop}
\label{new-properties} Let $H$ be a cocommutative weak Hopf algebra,
$(A,\varphi_{A})$ a weak left $H$-module algebra and $\sigma\in
Reg_{\varphi_{A}}(H^2,A)$.

\begin{itemize}

\item[(i)] $\sigma\co \Omega_{H}^{2}=\sigma.$

\item[(ii)] $\sigma_{H}^{A}\co \Omega_{H}^{2}=\sigma_{H}^{A}.$

\item[(iii)] $(A\ot \Omega_{H}^{2})\co (\sigma_{H}^{A}\ot H)=(\sigma_{H}^{A}\ot H)\co
(H\ot \Omega_{H}^{2}).$

\item[(iv)] $\partial_{2,3}(\sigma)=(\sigma\ot \varepsilon_{H}) \co
(H\ot \Omega_{H}^{2}).$

\end{itemize}

\end{prop}

{\em Proof}: To prove (i) first we show that $u_{2}\co
\Omega_{H}^{2}=u_{2}$. Indeed: By (\ref{omega-mu}) we have
$$u_{2}\co \Omega_{H}^{2}= \varphi_{A}\co ((\mu_{H}\co
\Omega_{H}^{2})\ot \eta_{A})=\varphi_{A}\co (\mu_{H}\ot
\eta_{A})=u_{2}.$$ Then, (i) holds because, by
(\ref{firts-cocom-1}), we obtain:
$$\sigma= \sigma\wedge \sigma^{-1}\wedge \sigma=\mu_{A}\co (u_{2}\ot \sigma)\co \delta_{H^{2}}
=\mu_{A}\co ((u_{2}\co \Omega_{H}^{2})\ot \sigma)\co
\delta_{H^{2}}$$
$$=\mu_{A}\co (u_{2}\ot \sigma)\co \delta_{H^{2}}\co \Omega_{H}^{2}=
(\sigma\wedge \sigma^{-1}\wedge \sigma)\co \Omega_{H}^{2}=\sigma\co
\Omega_{H}^{2}.$$

By (\ref{firts-cocom-1}) and the properties of (i) we have:
$$\sigma_{H}^{A}\co \Omega_{H}^{2}=((\sigma\co \Omega_{H}^{2})\ot
\mu_{H})\co \delta_{H^{2}}=\sigma_{H}^{A}.$$

Then (ii) holds.

Using that $\Omega_{H}^{2}$ is a morphism of left $H$-comodules and
$H$-modules we obtain (iii). Finally, (iv) is a consequence of
(\ref{omega-L-1}).

\begin{prop}
Let $H$ be a cocommutative weak Hopf algebra, $(A,\varphi_{A})$ a
weak left $H$-module algebra and $\sigma\in
Reg_{\varphi_{A}}(H^2,A)$. Then $\sigma$ satisfies 2-cocycle
condition if and only if the equality
\begin{equation}
\label{cocycle-equivalent} \mu_{A}\co (A\ot \sigma)\co
(\sigma_{H}^{A}\ot H)= \mu_{A}\co (A\ot \sigma)\co (\psi_{H}^{A}\ot
H)\co (H\ot \sigma_{H}^{A})
\end{equation}
holds.
\end{prop}

{\em Proof}: The proof follows from the following facts: First note

\begin{itemize}

\item[ ]$\hspace{0.38cm}\partial_{2,3}(\sigma)\wedge \partial_{2,1}(\sigma)$

\item [ ]$=\mu_{A}\co (((\sigma\ot \varepsilon_{H})\co
(H\ot \Omega_{H}^{2}))\ot (\sigma\co (\mu_{H}\ot H)))\co
\delta_{H^{3}}$

\item [ ]$=\mu_{A}\co (A\ot \sigma )\co (\sigma_{H}^{A}\ot H)\co (H\ot \Omega_{H}^{2})$

\item [ ]$=\mu_{A}\co (A\ot (\sigma\co \Omega_{H}^{2}) )\co (\sigma_{H}^{A}\ot H)$

\item [ ]$=\mu_{A}\co (A\ot \sigma )\co (\sigma_{H}^{A}\ot H)$

\end{itemize}

where the first equality follows by (iv) of Proposition
\ref{new-properties}, the second one by the properties of
$\varepsilon_{H}$, the third one by (iii) of Proposition
\ref{new-properties} and,  the last one by (i) of Proposition
\ref{new-properties}.

On  the other hand, by the naturality of $c$ we obtain that
$$\partial_{2,0}(\sigma)\wedge \partial_{2,2}(\sigma)=\mu_{A}\co (A\ot \sigma)\co (\psi_{H}^{A}\ot
H)\co (H\ot \sigma_{H}^{A})$$ and this finish the proof.

\begin{rem}
\label{twisted+cocycle+rem} {\rm Note that, if $(A, \varphi_{A})$ is
a commutative left $H$-module algebra, the 2-cocycle condition means
that $\sigma\in Ker(D^{2}_{\varphi_{A}})$.

Also, by the cocommutativity of $H$, we have
\begin{equation}
\label{sigma-commutative} \sigma_{H}^{A}=c_{A,H}\co \tau_{H}^{A}
\end{equation}

where $\tau_{H}^{A}=(\mu_{H}\ot \sigma)\co \delta_{H^{2}}$.
Therefore, if $(A, \varphi_{A})$ is a commutative left $H$-module
algebra the twisted condition holds for all $\sigma\in
Reg_{\varphi_{A}}(H^2,A)$.

}
\end{rem}

\begin{teo}
\label{prop-twisted} Let $H$ be a cocommutative weak Hopf algebra,
$(A,\varphi_{A})$  a weak left $H$-module algebra and $\sigma\in
Reg_{\varphi_{A}}(H^2,A)$.  The morphism $\sigma$ satisfies the
twisted condition (\ref{twisted-sigma}) if and only if ${\Bbb
A}_{H}$ satisfies the twisted condition (\ref{twis-wcp}).
\end{teo}

{\em Proof}: If ${\Bbb A}_{H}$ satisfies the twisted condition
(\ref{twis-wcp}), composing with $A\otimes \varepsilon_{H}$ and
using (ii) of Proposition \ref{sigma-prop} we obtain that $\sigma$
satisfies the twisted condition (\ref{twisted-sigma}). Conversely,
assume that $\sigma$ satisfies the twisted condition
(\ref{twisted-sigma}). Then

\begin{itemize}

\item[ ]$\hspace{0.38cm}(\mu_A\ot H)\co (A\ot \sigma_{H}^{A})\co (\psi_{H}^{A}\ot H)\co
(H\ot \psi_{H}^{A}) $

\item [ ]$= ((\mu_{A}\co ((\varphi_{A}\co (H\ot \varphi_{A}))\ot \sigma))\ot \mu_{H})\co
(H\ot ((H\ot A\ot H\ot c_{H,H})\co (H\ot A\ot c_{H,H}\ot H)\co $
\item[ ]$\hspace{0.38cm}  (H\ot c_{H,A}\ot H\ot H)\co (c_{H,H}\ot c_{H,A}\ot H)\co (H\ot
c_{H,H}\ot c_{H,A}))\ot H)\co (H\ot \delta_{H}\ot \delta_{H}\ot
c_{H,A})\co$
\item[ ]$\hspace{0.38cm}  (\delta_{H}\ot \delta_{H}\ot A)$

\item [ ]$=((\mu_{A}\co ((\varphi_{A}\co (H\ot \varphi_{A}))\ot A))\ot \mu_{H})\co
(H\ot H\ot A\ot c_{H,A}\ot H)\co (H\ot H\ot c_{H,A}\ot A\ot H)\co $
\item[ ]$\hspace{0.38cm} (H\ot c_{H,H}\ot c_{A,A}\ot H)\co
(c_{H,H}\ot H\ot \sigma\ot A\ot H)\co (H\ot \delta_{H^{2}}\ot
c_{H,A})\co (\delta_{H}\ot \delta_{H}\ot A)$

\item [ ]$=(A\ot \mu_{H})\co (c_{H,A}\ot H)\co   $
\item[ ]$\hspace{0.38cm} (H\ot (\mu_{A}\co ((\varphi_{A}\co (H\ot \varphi_{A}))\ot A)\co (H\ot H\ot c_{A,A})\co
(((H\ot H\ot \sigma)\co \delta_{H^{2}}) \ot A))\ot H)\co $
\item[ ]$\hspace{0.38cm}(H\ot H\ot H\ot c_{H,A})\co (\delta_{H}\ot
\delta_{H}\ot A)$

\item [ ]$=(A\ot \mu_{H})\co (c_{H,A}\ot H)\co    (H\ot (\mu_{A}\co (A\ot
\varphi_{A})\co (\sigma_{H}^{A}\ot A))\ot H)\co (H\ot H\ot H\ot
c_{H,A})\co $
\item[ ]$\hspace{0.38cm}(\delta_{H}\ot \delta_{H}\ot A)$

\item [ ]$= (\mu_{A}\ot H)\co (A\ot \varphi_{A}\ot \mu_{H})\co
(\sigma\ot \mu_{H}\ot c_{H,A}\ot H)\co (H\ot c_{H,H}\ot c_{H,H}\ot
c_{H,A})\co$
\item[ ]$\hspace{0.38cm}(H\ot H\ot  c_{H,H}\ot H\ot H\ot A)\co (
((H\ot \delta_{H})\co \delta_{H})\ot ((H\ot \delta_{H})\co
\delta_{H})\ot A)$

\item [ ]$= (\mu_{A}\ot H)\co (A\ot \varphi_{A}\ot H)\co (A\ot H\ot c_{H,A})\co
(\sigma\ot ((\mu_{H}\ot \mu_{H})\co \delta_{H^{2}})\ot A)\co
(\delta_{H^{2}}\ot A)$

\item[ ]$=(\mu_A\ot H)\co (A\ot \psi_{H}^{A})\co (\sigma_{H}^{A}\ot A).$

\end{itemize}

The first and the fifth equalities follow by the naturality of $c$,
the cocommutativity of $H$ and the coassociativity of $\delta_{H}$,
the second one by the cocommutativity of $H$ and the coassociativity
of $\delta_{H}$, the third and the sixth ones by the the naturality
of $c$, the fourth one by the twisted condition for $\sigma$ and the
last one by  (a1) of Definition (\ref{wha}).

Therefore, ${\Bbb A}_{H}$ satisfies the twisted condition
(\ref{twis-wcp}).

\begin{teo}
\label{prop-cocycle} Let $H$ be a cocommutative weak Hopf algebra,
$(A,\varphi_{A})$  a weak left $H$-module algebra and $\sigma\in
Reg_{\varphi_{A}}(H^2,A)$.  The morphism $\sigma$ satisfies the
2-cocycle condition (\ref{cocycle-equivalent}) if and only if ${\Bbb
A}_{H}$ satisfies the cocycle condition (\ref{cocy2-wcp}).
\end{teo}

{\em Proof}: If ${\Bbb A}_{H}$ satisfies the cocycle condition
(\ref{cocy2-wcp}), composing with $A\otimes \varepsilon_{H}$ and
using (ii) of Proposition \ref{sigma-prop} we obtain that $\sigma$
satisfies the 2-cocycle condition (\ref{cocycle-equivalent}).
Conversely, assume that $\sigma$ satisfies the 2-cocycle condition
(\ref{2-cocycle-sigma}). Then

\begin{itemize}

\item[ ]$\hspace{0.38cm}(\mu_A\ot H)\co (A\ot \sigma_{H}^{A})\co (\psi_{H}^{A}\ot H)\co
(H\ot\sigma_{H}^{A}) $

\item [ ]$= (\mu_{A}\ot H)\co (A\ot \sigma\ot \mu_{H})\co (\psi_{H}^{A}\ot c_{H,H}\ot H)
\co (H\ot c_{H,A}\ot H\ot H)\co (\delta_{H}\ot ((A\ot \delta_{H})\co
\sigma_{H}^{A})) $

\item [ ]$=(\mu_{A}\ot H)\co (A\ot \sigma\ot \mu_{H})\co (\psi_{H}^{A}\ot c_{H,H}\ot H)
\co (H\ot c_{H,A}\ot H\ot H)\co (\delta_{H}\ot ((\sigma_{H}^{A}\ot
\mu_{H})\co \delta_{H^{2}}))$

\item [ ]$=((\mu_{A}\co (A\ot \sigma)\co (\psi_{H}^{A}\ot H)\co (H\ot \sigma_{H}^{A}))
\ot (\mu_{H}\co (\mu_{H}\ot H)))\co  \delta_{H^{3}}$

\item [ ]$=((\mu_{A}\co (A\ot \sigma)\co (\sigma_{H}^{A}\ot H))\ot
(\mu_{H}\co (H\ot \mu_{H})))\co \delta_{H^{3}}$

\item [ ]$= (\mu_{A}\ot H)\co (A\ot \sigma\ot \mu_{H})\co (A\ot H\ot c_{H,H}\ot H)\co (((\sigma_{H}^{A}\ot
\mu_{H})\co \delta_{H^{2}})\ot \delta_{H})$

\item[ ]$=(\mu_A\ot H)\co (A\ot \sigma_{H}^{A}) \co (\sigma_{H}^{A}\ot H).$

\end{itemize}

The first equality follows by the naturality of $c$ and the
coassociativity of $\delta_{H}$, the second and the  sixth ones by
Proposition \ref{sigma-prop1}, the third and the fifth ones by the
naturality of $c$ and the associativity of $\mu_{H}$, the fourth one
by the 2-cocycle condition (\ref{cocycle-equivalent}).

\begin{rem}
{\rm By Theorems \ref{prop-twisted} and \ref{prop-cocycle} and
applying the general theory of weak crossed products, we have the
following: If $\sigma\in Reg_{\varphi_{A}}(H^2,A)$ satisfies the
twisted condition (\ref{twisted-sigma}) (equi\-va\-len\-tly
(\ref{cocycle-equivalent})) and the 2-cocycle condition
(\ref{2-cocycle-sigma}), the  quadruple ${\Bbb A}_{H}$ defined in
Remark \ref{quadruple-1} satisfies the twisted and the cocycle
conditions (\ref{twis-wcp}), (\ref{cocy2-wcp}) and therefore the
 induced product  is associative. Conversely, by Theorem
\ref{thm1-wcp}, we obtain that, if the product induced by the
quadruple ${\Bbb A}_{H}$ defined in Remark \ref{quadruple-1} is
associative, ${\Bbb A}_{H}$ satisfies the twisted and the cocycle
condition and, by Theorems \ref{prop-twisted} and
\ref{prop-cocycle}, $\sigma$ satisfies the twisted condition
(\ref{twisted-sigma}) and the 2-cocycle condition
(\ref{2-cocycle-sigma}) (equivalently (\ref{cocycle-equivalent})).

}
\end{rem}

\begin{defin}
 {\rm Let $H$ be a  cocommutative weak Hopf
algebra, $(A,\varphi_{A})$  a weak left $H$-module algebra and
$\sigma\in Reg_{\varphi_{A}}(H^2,A)$.  We say that $\sigma$
satisfies the normal condition if
\begin{equation}
\label{normal-sigma} \sigma\co (\eta_{H}\ot H)=\sigma\co (H\ot
\eta_{H})=u_{1},
\end{equation}
i.e. $\sigma\in Reg_{\varphi_{A}}^+(H^2,A)$.
 }
\end{defin}

\begin{teo}
\label{norma-sigma-prop3} Let $H$ be a cocommutative weak Hopf
algebra, $(A,\varphi_{A})$  a weak left $H$-module algebra and
$\sigma\in Reg_{\varphi_{A}}(H^2,A)$. Let  ${\Bbb A}_{H}$ be the
quadruple defined in Remark \ref{quadruple-1} and assume that ${\Bbb
A}_{H}$ satisfies the twisted and the cocycle conditions
(\ref{twis-wcp}) and (\ref{cocy2-wcp}). Then, $\nu=\nabla_{A\ot
H}\co (\eta_{A}\ot\eta_{H})$ is a preunit for the weak crossed
product associated to ${\Bbb A}_{H}$ if and only if
\begin{equation}
\label{sigma-preunit1} \sigma_{H}^{A}\co (\eta_{H}\ot
H)=\sigma_{H}^{A}\co (H\ot \eta_{H})=\nabla_{A\ot H}\co (\eta_{A}\ot
H).
\end{equation}

\end{teo}

{\em Proof}: By Theorem \ref{thm1-wcp}, to prove the result, we only
need to show that (\ref{pre1-wcp}), (\ref{pre2-wcp}) and
(\ref{pre3-wcp}) hold for $\nu=\nabla_{A\ot H}\co
(\eta_{A}\ot\eta_{H})$ if and only if
$$\sigma_{H}^{A}\co (\eta_{H}\ot
H)=\sigma_{H}^{A}\co (H\ot \eta_{H})=\nabla_{A\ot H}\co (\eta_{A}\ot
H).$$

Indeed, $\nu$ satisfies (\ref{pre1-wcp}) if and only if
$\sigma_{H}^{A}\co (H\ot \eta_{H})=\nabla_{A\ot H}\co (\eta_{A}\ot
H)$ because:
\begin{itemize}

\item[ ]$\hspace{0.38cm}(\mu_A\otimes H)\circ (A\otimes \sigma_{H}^{A})\circ
    (\psi_{H}^{A}\otimes H)\circ (H\otimes \nu) $

\item [ ]$=(\mu_A\otimes H)\circ (A\otimes \sigma_{H}^{A})\circ
    (\psi_{H}^{A}\otimes H)\circ (H\ot (\psi_{H}^{A}\co (\eta_{H}\ot \eta_{A})))  $

\item [ ]$=\nabla_{A\ot H}\co \sigma_{H}^{A}\co (H\ot \eta_{H})$

\item [ ]$= \sigma_{H}^{A}\co (H\ot \eta_{H}).$

\end{itemize}

The first equality follows by the definition of $\nabla_{A\ot H}$,
the second one by the twisted condition and the last one by (ii) of
Proposition \ref{sigma-prop}.

Also, $\nu$ satisfies (\ref{pre2-wcp}) if and only if
$\sigma_{H}^{A}\co ( \eta_{H}\ot H)=\nabla_{A\ot H}\co (\eta_{A}\ot
H)$ because by (\ref{aw1}) we have
$$(\mu_A\otimes H)\circ (A\otimes \sigma_{H}^{A})\circ
    (\nu\otimes H) = \sigma_{H}^{A}\co (\eta_{H}\ot H).$$

Finally, (\ref{pre3-wcp})  is always true because, by
(\ref{fi-nab}), we obtain
$$\hspace{0.38cm}(\mu_A\otimes H)\circ (A\otimes \psi_{H}^{A})\circ
(\nu\otimes A)=\psi_{H}^{A}\co (\eta_{H}\ot A)=\beta_{\nu}. $$

\begin{cor}
\label{norma-sigma-prop4} Let $H$ be a cocommutative weak Hopf
algebra, $(A,\varphi_{A})$  a weak left $H$-module algebra and
$\sigma\in Reg_{\varphi_{A}}(H^2,A)$. Let  ${\Bbb A}_{H}$ be the
quadruple defined in Remark \ref{quadruple-1} and assume that ${\Bbb
A}_{H}$ satisfies the twisted and the cocycle conditions
(\ref{twis-wcp}) and (\ref{cocy2-wcp}). Then, $\nu=\nabla_{A\ot
H}\co (\eta_{A}\ot\eta_{H})$ is a preunit for the weak crossed
product associated to ${\Bbb A}_{H}$ if and only if $\sigma$
satisfies the normal condition (\ref{normal-sigma}).
\end{cor}

{\em Proof}: If $\nu=\nabla_{A\ot H}\co (\eta_{A}\ot\eta_{H})$ is a
preunit for the weak crossed product associated to ${\Bbb A}_{H}$,
by Theorem \ref{norma-sigma-prop3} we have (\ref{sigma-preunit1}).
Then, composing whit $(A\ot \varepsilon_{H})$ and using (ii) of
Proposition \ref{sigma-prop}, we obtain (\ref{normal-sigma}).
Conversely, if (\ref{normal-sigma}) holds, we have:

\begin{itemize}

\item[ ]$\hspace{0.38cm}\sigma_{H}^{A}\co (\eta_{H}\ot H)$

\item [ ]$=((\sigma\co c_{H,H})\ot \mu_{H})\co (H\ot (\delta_{H}\co \eta_{H})\ot H)\co \delta_{H} $

\item [ ]$= ((\sigma\co c_{H,H})\ot H)\co (H\ot ((\overline{\Pi}_{H}^{L}\ot H)\co
\delta_{H}))\co \delta_{H}$

\item [ ]$=((\sigma\co c_{H,H}\co(H\ot \overline{\Pi}_{H}^{L})\co \delta_{H})\ot H)\co \delta_{H} $

\item [ ]$= (u_{1}\ot H)\co \delta_{H}$

\item [ ]$=\nabla_{A\otimes H}\co (\eta_{A}\ot H). $

\end{itemize}

The first equality follows by the naturality of $c$, the second one
by (\ref{delta-pi31}), the fourth one by the coassociativity of
$\delta_{H}$ and the fourth one by (i) of Proposition
\ref{new-properties-sigma}. The last one follows by definition.

On the other hand

\begin{itemize}

\item[ ]$\hspace{0.38cm}\sigma_{H}^{A}\co (H\ot \eta_{H})$

\item [ ]$=(\sigma\ot H)\co (H\ot ((\Pi_{H}^{R}\ot H)\co \delta_{H}))\co \delta_{H} $

\item [ ]$=((\sigma\co (H\ot \Pi_{H}^{R})\co \delta_{H})\ot H)\co \delta_{H} $

\item [ ]$= (u_{1}\ot H)\co \delta_{H}$

\item [ ]$=\nabla_{A\otimes H}\co (\eta_{A}\ot H). $

\end{itemize}

The first equality follows by (\ref{delta-pi4}), the second one by
the coassociativity of $\delta_{H}$,  the third one by (ii) of
Proposition \ref{new-properties-sigma}, and the last one  by
definition.

\begin{cor}
\label{crossed-product1}
 Let $H$ be a cocommutative weak Hopf
algebra, $(A,\varphi_{A})$  a weak left $H$-module algebra and
$\sigma\in Reg_{\varphi_{A}}(H^2,A)$. Let  ${\Bbb A}_{H}$ be the
quadruple defined in Remark \ref{quadruple-1} and $\mu_{A\otimes H}$
the associated product defined in (\ref{prod-todo-wcp}). Then the
following statements are equivalent:
\begin{itemize}

\item[(i)] The product $\mu_{A\otimes H}$ is associative with preunit
$\nu=\nabla_{A\otimes H}\co (\eta_{A}\ot \eta_{H})$ and normalized
with respect to $\nabla_{A\otimes H}.$

\item[(ii)] The morphism $\sigma$ satisfies the twisted condition (\ref{twisted-sigma}),
the 2-cocycle condition (\ref{2-cocycle-sigma}) (equivalently
(\ref{cocycle-equivalent})) and  the normal condition
(\ref{normal-sigma}).

\end{itemize}

\end{cor}

{\em Proof}: The proof is an easy consequence of Theorems
\ref{thm1-wcp}, \ref{prop-twisted}, \ref{prop-cocycle} and Corollary
\ref{norma-sigma-prop4}.

\begin{notac}
{\rm Let $H$ be a cocommutative weak Hopf algebra and
$(A,\varphi_{A})$  a weak left $H$-module algebra. From now on we
will denote by $A\ot_{\tau} H=(A\ot H, \mu_{A\ot_{\tau}H})$ the weak
crossed product, with preunit $\nu=\nabla_{A\ot H}\co
(\eta_{A}\ot\eta_{H})$, defined by  $\tau\in
Reg_{\varphi_{A}}(H^2,A)$ when it satisfies the twisted condition,
the 2-cocycle condition and  the normal condition. The associated
algebra will be denoted by
$$A\times_{\tau} H=(A\times H,
\eta_{A\times_{\tau}H}, \mu_{A\times_{\tau}H}).$$

 Finally, the quadruple ${\Bbb A}_{H}$ defined
in Remark \ref{quadruple-1} will be denoted by ${\Bbb A}_{H,\tau}$
and $\sigma_{H}^{A}$ by $\sigma_{H,\tau}^{A}$.

}
\end{notac}

\begin{rem}
{\rm Let $H$ be a cocommutative weak Hopf algebra and
$(A,\varphi_{A})$  a weak left $H$-module algebra. Let  $\sigma\in
Reg_{\varphi_{A}}(H^2,A)$ be  a morphism satisfying the twisted
condition (\ref{twisted-sigma}), the 2-cocycle condition
(\ref{2-cocycle-sigma}) and  the normal condition
(\ref{normal-sigma}). Then, the weak crossed product $A\ot_{\sigma}
H=(A\ot H, \mu_{A\ot_{\sigma}H})$  with preunit $\nu=\nabla_{A\ot
H}\co (\eta_{A}\ot\eta_{H})$ defined previously is a particular
instance of the weak crossed products introduced in
\cite{mra-preunit}. Also is a particular case  of the ones used in
\cite{ana1} where these crossed structures were studied in a
category of modules over a commutative ring without requiring
cocommutativity of $H$ and using weak measurings (see Definition 3.2
of \cite{ana1}). To prove this assertion we will show that the
conditions presented in Lemma 3.8 and Theorem 3.9 of \cite{ana1} are
completely fulfilled. First, note that, if $(A,\varphi_{A})$  a weak
left $H$-module algebra, we have that $\varphi_{A}$ is a weak
measuring. The idempotent morphism $\Omega_{A\ot H}$ related with
the preunit $\nu$ is the morphism $\nabla_{A\ot H}$ because, by
(\ref{aw1}) and (\ref{sigma-preunit1}), we have
$$\Omega_{A\ot H}=\mu_{A\ot_{\sigma}H}\co (A\ot H\ot \nu)=\mu_{A\ot_{\sigma}H}\co (A\ot H\ot \eta_A\ot \eta_H)=(\mu_{A}\ot H)\co (A\ot \sigma_{H}^{A})\co (\nabla_{A\ot H}\ot \eta_{H})$$
$$=(\mu_{A}\ot H)\co (A\ot ( \nabla_{A\ot H}\co (\eta_{A}\ot H)))=\nabla_{A\ot H}.$$
Moreover, in the category of modules over and associative
commutative unital ring, the normalized condition implies that
$Im(\mu_{A\ot_{\sigma}H})\subset Im(\nabla_{A\ot H})$.

On the other hand, the left action defined in Lemma 3.8 of
\cite{ana1} is $\varphi_{A}$. Indeed:

\begin{itemize}

\item[ ]$\hspace{0.38cm}(A\ot \varepsilon_{H})\co\mu_{A\ot_{\sigma}H}\co (\eta_{A}\ot H\ot
((\mu_{A}\ot H)\co (A\ot \nu)))$

\item [ ]$=(\mu_A\ot \varepsilon_{H})\co (A\ot \sigma_{H}^{A})\co (\psi_{H}^{A}\ot H)\co (H\ot (\nabla_{A\ot H}\co (A\ot \eta_{H})))$

\item [ ]$=(\mu_A\ot \varepsilon_{H})\co (A\ot \sigma_{H}^{A})\co (\psi_{H}^{A}\ot H)\co (H\ot A\ot \eta_{H}) $

\item [ ]$=\mu_{A}\co (A\ot \sigma)\co  (\psi_{H}^{A}\ot \eta_{H})$

\item [ ]$=\mu_{A}\co (A\ot u_{1})\co  \psi_{H}^{A}$

\item [ ]$=\varphi_{A}. $

\end{itemize}

where the first equality follows by the unit properties, the second
one by  (\ref{c11}), the third one by (iii) of Proposition
\ref{sigma-prop}, the fourth one by (\ref{normal-sigma}) and finally
the last one (b2) of Definition \ref{weak-H-mod}.

Also, the morphism defined in Lemma 3.8 of \cite{ana1} is $\sigma$
because, by (\ref{aw1}) and (iii) of Proposition (\ref{sigma-prop}),
we have
$$(A\ot \varepsilon_{H})\co\mu_{A\ot_{\sigma}H}\co (\eta_{A}\ot H\ot \eta_{A}\ot H)=
(\mu_A\ot \varepsilon_{H})\co (A\ot \sigma_{H}^{A})\co
((\nabla_{A\ot H}\co (\eta_{A}\ot H))\ot H)$$
$$=\mu_{A}\co (\eta_{A}\ot \sigma)=\sigma.$$

Then, the equalities (a) and (b) of Lemma 3.8 of \cite{ana1} hold
because the first one is the definition of $\psi_{H}^{A}$ and the
second one is a consequence of (\ref{aw1}) and the definition of
$\sigma_{H}^{A}$.

Therefore, we have that $\mu_{A\ot_{\sigma}H}$ satisfies that
$$\rho_{A\ot H}\co \mu_{A\ot_{\sigma}H}=(\mu_{A\ot_{\sigma}H}\ot H)\co \rho_{A\ot H\ot A\ot H}$$
where $\rho_{A\ot H}=A\ot \delta_{H}$ and
$$\rho_{A\ot H\ot A\ot H}=(A\ot H\ot A\ot H\ot \mu_{H})\co (A\ot H\ot c_{H,A\ot H}\ot H)\co (\rho_{A\ot H}\ot \rho_{A\ot H}).$$
Although that $\rho_{A\ot H\ot A\ot H}$ it is not counital, we say
that $\mu_{A\ot_{\sigma}H}$ is $H$-colinear as in Lemma 3.8 of
\cite{ana1}. Then we obtain that $\sigma$ satisfies the equality (1)
of \cite{ana1}, that is:
$$\sigma\co ((\mu_{H}\co (H\ot \Pi_{H}^{R}))\ot H)=\sigma\co (H\ot (\mu_{H}\co (\Pi_{H}^{R}\ot H))).$$

Finally, for the preunit $\nu=\nabla_{A\ot H}\co
(\eta_{A}\ot\eta_{H})$, by the equalities (\ref{nabla-nabla}) and
(\ref{delta-pi2}),
$$(A\ot \delta_{H})\co \nu=(A\ot ((H\ot \Pi_{H}^L)\co\delta_{H}))\co \nu$$
holds (i.e., the equality (4) of \cite{ana1}  is true in our
setting).

}
\end{rem}

\section{Equivalent weak crossed products and
$H^{2}_{\varphi_{A}}(H,A)$}

The aim of this section is to give necessary and sufficient
conditions for two weak crossed products $A\ot_{\alpha} H$,
$A\ot_{\beta} H$ to be equivalent in the cocommutative setting. To
define a good notion of equivalence we need the definition of right
$H$-comodule algebra for a weak Hopf algebra $H$.

\begin{defin}
\label{pre} {\rm Let $H$ be a weak bialgebra and $(B, \rho_{B})$ an
algebra which is also a right $H$-comodule such that
\begin{equation}
\label{comod-alg} \mu_{B\otimes H}\circ (\rho_{B}\otimes
\rho_{B})=\rho_{B}\circ \mu_{B}.
\end{equation}

 The object $(B, \rho_{B})$ is called a right
$H$-comodule algebra if one of the following equivalent conditions
holds:

\begin{itemize}

\item[(d1)]$(\rho_{B}\otimes H)\circ \rho_{B}\circ
\eta_{B}=(B\otimes (\mu_{H}\circ c_{H,H})\otimes H)\circ
((\rho_{B}\circ \eta_{B})\otimes (\delta_{H}\circ \eta_{H})),$

\item[(d2)]$(\rho_{B}\otimes H)\circ \rho_{B}\circ
\eta_{B}=(B\otimes \mu_{H}\otimes H)\circ ((\rho_{B}\circ
\eta_{B})\otimes (\delta_{H}\circ \eta_{H})),$

\item[(d3)]$(B\otimes \overline{\Pi}_{H}^{R})\circ
\rho_{B}=(\mu_{B}\otimes H)\circ (B\otimes (\rho_{B}\circ
\eta_{B})),$

 \item[(d4)]$(B\otimes \Pi_{H}^{L})\circ
\rho_{B}=((\mu_{B}\circ c_{B,B})\otimes H)\circ (B\otimes
(\rho_{B}\circ \eta_{B})),$

\item[(d5)]$(B\otimes
\overline{\Pi}_{H}^{R})\circ \rho_{B}\circ \eta_{B}=\rho_{B}\circ
\eta_{B},$

\item[(d6)]$(B\otimes \Pi_{H}^{L})\circ \rho_{B}\circ
\eta_{B}=\rho_{B}\circ \eta_{B}.$

\end{itemize}

}
\end{defin}

\begin{prop}
\label{H-comod-alg} Let $H$ be a cocommutative weak Hopf algebra,
$(A,\varphi_{A})$  a weak left $H$-module algebra and $\alpha\in
Reg_{\varphi_{A}}^{+}(H^2,A)$ such that satisfies the twisted
condition (\ref{twisted-sigma}) and the 2-cocycle condition
(\ref{2-cocycle-sigma}) (equivalently (\ref{cocycle-equivalent})).
Then, the algebra $A\times_{\alpha} H=(A\times H,
\eta_{A\times_{\alpha}H}, \mu_{A\times_{\alpha}H})$ is a right
$H$-comodule algebra for the coaction
$$\rho_{A\times_{\alpha}H}=(p_{A\ot H}\ot H)\co (A\ot \delta_{H})\co i_{A\ot H}.$$

\end{prop}

{\em Proof}: First note that $(A\times_{\alpha}H,
\rho_{A\times_{\alpha}H})$ is a right $H$-comodule because

$$(A\times H\ot \varepsilon_{H})\co \rho_{A\times_{\alpha}H}=p_{A\ot H}\co i_{A\ot H}
=id_{A\times H}$$  and, by (\ref{nabla-delta}) and the
coassociativity of $\delta_{H}$,
$$(\rho_{A\times_{\alpha}H}\ot H)\co \rho_{A\times_{\alpha}H}=(p_{A\ot H}\ot H\ot H)
\co (A\ot ((\delta_{H}\ot H)\co \delta_{H}))\co i_{A\ot H}=(A\times
H\ot \delta_{H}) \co \rho_{A\times_{\alpha}H}.$$

On the other hand,

\begin{itemize}

\item[ ]$\hspace{0.38cm} \mu_{(A\times_{\alpha}H)\ot H}\co (\rho_{A\times_{\alpha}H}\ot  \rho_{A\times_{\alpha}H}) $

\item [ ]$= (p_{A\times H}\ot H)\co (\mu_{A\otimes_{\alpha}H}\ot \mu_{H})\co (A\ot H\ot A\ot c_{H,H}\ot H)\co (A\ot H\ot c_{H,A}\ot H\ot H)\co $
\item[ ]$\hspace{0.38cm} (((A\ot \delta_{H})\co i_{A\ot H})\ot ((A\ot \delta_{H})\co i_{A\ot H}))$

\item [ ]$= (p_{A\times H}\ot H)\co (\mu_{A}\ot H\ot H)\co (\mu_{A}\ot ((\sigma_{H}^{A}\ot
H\ot H)\co \delta_{H^{2}}))\co (A\ot \psi_{H}^{A}\ot H)\co (i_{A\ot
H}\ot i_{A\ot H}) $

\item [ ]$= (p_{A\times H}\ot H)\co (\mu_{A}\ot H\ot H)\co (\mu_{A}\ot ((A\ot \delta_{H})\co \sigma_{H}^{A}))\co (A\ot \psi_{H}^{A}\ot H)\co (i_{A\ot H}\ot i_{A\ot H}) $

\item [ ]$=  \rho_{A\times_{\alpha}H}\co \mu_{A\times_{\alpha}H} $

\end{itemize}

where the first equality follows by the normalized condition for
$\mu_{A\ot_{\alpha}H}$, the second one by the naturality of $c$ and
the coassociativity of $\delta_{H}$, the third one by
(\ref{delta-sigmaHA}) and the last one by (\ref{nabla-delta}).

Finally, by (\ref{nabla-delta}) and (\ref{delta-pi2}), we obtain
that
$$(A\times_{\alpha}H\ot \Pi_{H}^{L} )\co \rho_{A\times_{\alpha}H}\co
\eta_{A\times_{\alpha}H}=(p_{A\times H}\ot \Pi_{H}^{L} )\co
(\eta_{A}\ot (\delta_{H}\co \eta_{H}))=(p_{A\times H}\ot H)\co
(\eta_{A}\ot (\delta_{H}\co \eta_{H}))$$
$$=\rho_{A\times_{\alpha}H}\co \eta_{A\times_{\alpha}H}.$$

\begin{defin}
{\rm Let $H$ be a cocommutative weak Hopf algebra, $(A,\varphi_{A})$
a weak left $H$-module algebra and $\alpha, \beta\in
Reg_{\varphi_{A}}^{+}(H^2,A)$ such that satisfy the twisted
condition (\ref{twisted-sigma}) and the 2-cocycle condition
(\ref{2-cocycle-sigma}) (equivalently (\ref{cocycle-equivalent})).
Let $A\ot_{\alpha} H$, $A\ot_{\beta} H$ be the weak crossed products
associated to $\alpha$ and $\beta$. We say that $A\ot_{\alpha} H$,
$A\ot_{\beta} H$ are equivalent if there is an isomorphism of left
$A$-modules and right $H$-comodule algebras
$\omega_{\alpha,\beta}:A\times_{\alpha} H\rightarrow A\times_{\beta}
H$. }
\end{defin}

\begin{rem}
\label{relevant-remark} {\rm Let  $H$ be a  weak Hopf algebra,
$(A,\varphi_{A})$ a weak left $H$-module algebra. Let $\Gamma:A\ot
H\rightarrow A\ot H$ be a morphism of left $A$-modules and right
$H$-comodules for the regular action $\varphi_{A\otimes
H}=\mu_{A}\ot H$ and coaction $\rho_{A\ot H}=A\ot \delta_{H}$. Then
\begin{equation}
\label{aux-import} \Gamma\co (\eta_{A}\ot H)=(A\ot
\varepsilon_{H}\ot H)\co \rho_{A\ot H}\co \Gamma\co (\eta_{A}\ot
H)=(f_{\Gamma}\ot H)\co \delta_{H} \end{equation} where
$f_{\Gamma}=(A\ot \varepsilon_{H})\co \Gamma\co (\eta_{A}\ot H)$. As
a consequence:
\begin{equation}
\label{f-gamma} \Gamma=(\mu_{A}\ot H)\co (A\ot (\Gamma\co
(\eta_{A}\ot H)))= ((\mu_{A}\co (A\ot f_{\Gamma}))\ot H)\co (A\ot
\delta_{H}).
\end{equation}

If $f:H\rightarrow A$ is a morphism and we define
$$\Gamma_{f}:A\ot H\rightarrow A\ot H$$

 by $\Gamma_{f}=((\mu_{A}\co (A\ot f))\ot H)\co (A\ot
\delta_{H})$, it is clear that $\Gamma_{f}$ is a morphism of left
$A$-modules and right $H$-comodules such that $f_{\Gamma_{f}}=f$.
Also, $\Gamma_{f_{\Gamma}}=\Gamma$ and then there is a bijection
$$\Phi:\;_{A}Hom^{H}_{\mathcal C}(A\ot H, A\ot H)\rightarrow Hom_{\mathcal C}(H, A)$$
defined by $\Phi(\Gamma)=f_{\Gamma}$ which inverse
$\Phi^{-1}(f)=\Gamma_{f}$. Note that
$\Phi^{-1}(u_{1})=\Gamma_{u_{1}}=\nabla_{A\ot H}.$

Then, it is easy to show  that $\Gamma, \Gamma^{\prime} \in
\;_{A}Hom^{H}_{\mathcal C}(A\ot H, A\ot H)$ satisfy
\begin{itemize}

\item[(e1)] $\Gamma\co \Gamma^{\prime}=\Gamma^{\prime}\co \Gamma=\nabla_{A\ot H}.$

\item[(e2)] $\Gamma\co \Gamma^{\prime}\co \Gamma=\Gamma.$

\item[(e3)] $\Gamma^{\prime}\co \Gamma\co \Gamma^{\prime}=\Gamma^{\prime}.$

\end{itemize}
 if and only if for the morphism $f_{\Gamma}$ there exists a morphism $f_{\Gamma}^{-1}$ satisfying:
\begin{itemize}

\item[(i)] $f_{\Gamma}\wedge f_{\Gamma}^{-1}= f_{\Gamma}^{-1}\wedge
f_{\Gamma}=u_{1}.$

\item[(ii)] $f_{\Gamma}\wedge f_{\Gamma}^{-1}\wedge f_{\Gamma}=f_{\Gamma}.$

\item[(iii)] $f_{\Gamma}^{-1}\wedge f_{\Gamma}\wedge f_{\Gamma}^{-1}=f_{\Gamma}^{-1}.$

\end{itemize}
Indeed: If $\Gamma, \Gamma^{\prime} \in \;_{A}Hom^{H}_{\mathcal
C}(A\ot H, A\ot H)$ satisfies (e1)-(e3) define $f_{\Gamma}^{-1}$ by
$f_{\Gamma}^{-1}=f_{\Gamma^{\prime}}, $ and, conversely, if for
$f_{\Gamma}$ there exists a morphism $f_{\Gamma}^{-1}$ satisfying
(i)-(iii), define $\Gamma^{\prime}$ by
$\Gamma^{\prime}=\Gamma_{f_{\Gamma}^{-1}}.$

As a consequence, if $H$ is cocommutative, $\Gamma \in
\;_{A}Hom^{H}_{\mathcal C}(A\ot H, A\ot H)$ satisfies (e1)-(e3) if
and only if $\Phi(\Gamma)=f_{\Gamma}\in Reg_{\varphi_{A}}(H,A).$
Conversely, $f\in Reg_{\varphi_{A}}(H,A)$ if and only if
$\Phi^{-1}(f)=\Gamma_{f}$ satisfies (e1)-(e3). }
\end{rem}

\begin{teo}
\label{principal-1} Let $H$ be a cocommutative weak Hopf algebra,
$(A,\varphi_{A})$ a weak left $H$-module algebra and $\alpha,
\beta\in Reg_{\varphi_{A}}^{+}(H^2,A)$ such that satisfy the twisted
condition (\ref{twisted-sigma}) and the 2-cocycle condition
(\ref{2-cocycle-sigma}) (equivalently (\ref{cocycle-equivalent})).
The weak crossed products $A\ot_{\alpha} H$, $A\ot_{\beta} H$
associated to $\alpha$ and $\beta$ are equivalent if and only if
there exist  multiplicative and preunit preserving morphisms
$\Gamma, \Gamma^{\prime} \in \;_{A}Hom^{H}_{\mathcal C}(A\ot H, A\ot
H)$ satisfying (e1)-(e3).
\end{teo}
{\em Proof}: Assume that $A\ot_{\alpha} H$, $A\ot_{\beta} H$ are
equivalent. Thus there exists and isomorphism of left $A$-modules
and right $H$-comodule algebras
$\omega_{\alpha,\beta}:A\times_{\alpha} H\rightarrow A\times_{\beta}
H$. Define $\Gamma$ and $\Gamma^{\prime}$ by
$$\Gamma=i_{A\ot H}\co \omega_{\alpha,\beta}\co p_{A\ot H},\;\;\;\Gamma^{\prime}=
i_{A\ot H}\co \omega_{\alpha,\beta}^{-1}\co p_{A\ot H}.$$

Then,
$$\Gamma\co \Gamma^{\prime}=i_{A\ot H}\co \omega_{\alpha,\beta}\co p_{A\ot H}\co i_{A\ot H}\co
\omega_{\alpha,\beta}^{-1}\co p_{A\ot H}=i_{A\ot H}\co
\omega_{\alpha,\beta}\co \omega_{\alpha,\beta}^{-1}\co p_{A\ot
H}=\nabla_{A\ot H},$$ and
$$\Gamma^{\prime}\co \Gamma=i_{A\ot H}\co \omega_{\alpha,\beta}^{-1}\co p_{A\ot H}\co
 i_{A\ot H}\co \omega_{\alpha,\beta}\co p_{A\ot H}=i_{A\ot H}\co \omega_{\alpha,\beta}^{-1}
 \co \omega_{\alpha,\beta}\co p_{A\ot H}=\nabla_{A\ot H}.$$

Also,
$$\Gamma\co \Gamma^{\prime}\co \Gamma=\nabla_{A\ot H}\co \Gamma=\Gamma, \;\;\;
\Gamma^{\prime}\co \Gamma\co \Gamma^{\prime}=\nabla_{A\ot H}\co
\Gamma^{\prime}= \Gamma^{\prime}$$ and therefore (e1)-(e3) hold.

The morphism $\Gamma$ is multiplicative because
$\omega_{\alpha,\beta}$ is an algebra morphism:

\begin{itemize}

\item[ ]$\hspace{0.38cm} \mu_{A\ot_{\beta}H}\co (\Gamma\ot \Gamma)$

\item [ ]$= \mu_{A\ot_{\beta}H}\co (i_{A\ot H}\ot i_{A\ot H})\co (\omega_{\alpha,\beta}\ot
\omega_{\alpha,\beta})\co (p_{A\ot H}\ot p_{A\ot H})$

\item [ ]$= i_{A\ot H}\co \mu_{A\times_{\beta}H}\co (\omega_{\alpha,\beta}\ot
\omega_{\alpha,\beta})\co (p_{A\ot H}\ot p_{A\ot H}) $

\item [ ]$= i_{A\ot H}\co \omega_{\alpha,\beta}\co \mu_{A\times_{\alpha}H}
\co (p_{A\ot H}\ot p_{A\ot H})$

\item [ ]$= \Gamma\co  \mu_{A\ot_{\alpha}H} $

\end{itemize}

and in a similar way, using that $\omega_{\alpha,\beta}^{-1}$ is
multiplicative, it is possible to prove that $\Gamma^{\prime}$ is
multiplicative.

On the other hand, $\Gamma$ preserve the preunit because:
$$\Gamma\co \nu=i_{A\ot H}\co \omega_{\alpha,\beta}\co \eta_{A\times_{\alpha}H}=
i_{A\ot H}\co  \eta_{A\times_{\beta}H}=\nu.$$ By the same arguments
we obtain that $\Gamma^{\prime}\co \nu=\Gamma^{\prime}.$

Using (e1), (e2) and the left $A$-linearity of
$\omega_{\alpha,\beta}$  we have

\begin{itemize}

\item[ ]$\hspace{0.38cm} \varphi_{A\ot H}\co (A\ot \Gamma)$

\item [ ]$=\varphi_{A\ot H}\co (A\ot (\nabla_{A\ot H}\co\Gamma)) $

\item [ ]$=\nabla_{A\ot H}\co (\mu_{A}\ot H)\co (A\ot \Gamma) $

\item [ ]$= i_{A\ot H}\co \varphi_{A\times_{\beta}H}\co (A\ot \omega_{\alpha,\beta})
\co (A\ot p_{A\ot H})$

\item [ ]$= i_{A\ot H}\co \omega_{\alpha,\beta}\co \varphi_{A\times_{\alpha}H}\co
(A\ot p_{A\ot H})$

\item [ ]$= \Gamma \co (\mu_{A}\ot H)\co (A\ot \nabla_{A\ot H})$

\item [ ]$= \Gamma \co \nabla_{A\ot H}\co (\mu_{A}\ot H)$

\item [ ]$= \Gamma \co \varphi_{A\ot H}.$

\end{itemize}

Similarly, by (e1), (e3) and the left $A$-linearity of
$\omega_{\alpha,\beta}^{-1}$ we obtain that $\Gamma^{\prime}$ is a
morphism of left $A$-modules.

Finally, $\Gamma$ is a morphism of right $H$-comodules by
(\ref{nabla-delta}) and the  right $H$-comodule morphism property of
$\omega_{\alpha,\beta}$. Indeed:
\begin{itemize}

\item[ ]$\hspace{0.38cm} \rho_{A\ot H}\co \Gamma$

\item [ ]$= (i_{A\ot H}\ot H)\co \rho_{A\times_{\beta}H}\co \omega_{\alpha,\beta}
\co p_{A\ot H}$

\item [ ]$=((i_{A\ot H}\co\omega_{\alpha,\beta}) \ot H)\co \rho_{A\times_{\alpha}H}
\co  p_{A\ot H} $

\item [ ]$=(\Gamma\ot H)\co (A\ot \delta_{H})\co \nabla_{A\ot H} $

\item [ ]$= ((\Gamma\co \nabla_{A\ot H}) \ot H)\co (A\ot \delta_{H})$

\item [ ]$= (\Gamma \ot H)\co \rho_{A\ot H}$

\end{itemize}

By a similar calculus we obtain that $\Gamma^{\prime}$ is a morphism
of right $H$-comodules.

Conversely, assume that there exist  multiplicative and preunit
preserving morphisms
$$\Gamma, \Gamma^{\prime} \in \;_{A}Hom^{H}_{\mathcal C}(A\ot H, A\ot H)$$
satisfying (e1)-(e3) of the previous remark. Define
$$\omega_{\alpha,\beta}=p_{A\ot H}\co \Gamma\co i_{A\ot H},\;\;\;
\omega_{\alpha,\beta}^{-1}=p_{A\ot H}\co \Gamma^{\prime}\co i_{A\ot
H}.$$

Then, by (e1), (e2) and (e3), we have
$$\omega_{\alpha,\beta}^{-1}\co \omega_{\alpha,\beta}=p_{A\ot H}\co
\Gamma^{\prime}\co \nabla_{A\ot H}\co \Gamma\co i_{A\ot H}=p_{A\ot
H}\co \Gamma^{\prime}\co \Gamma\co i_{A\ot H}= p_{A\ot H}\co
\nabla_{A\ot H}\co i_{A\ot H}=id_{A\times H}$$ and
$$\omega_{\alpha,\beta}\co \omega_{\alpha,\beta}^{-1}=p_{A\ot H}\co
\Gamma\co \nabla_{A\ot H}\co \Gamma^{\prime}\co i_{A\ot H}= p_{A\ot
H}\co \Gamma\co \Gamma^{\prime}\co i_{A\ot H}=p_{A\ot H}\co
\nabla_{A\ot H}\co i_{A\ot H}=id_{A\times H}$$

which proves that $\omega_{\alpha,\beta}$ is an isomorphism.

Moreover, using that $\Gamma$ preserves the preunit
$\nu=\nabla_{A\ot H}\co (\eta_{A}\ot \eta_{H})$ we have
$$\omega_{\alpha,\beta}\co \eta_{A\times_{\alpha}H}=p_{A\ot H}\co
\Gamma\co\nu=p_{A\ot H}\co \nu=\eta_{A\times_{\beta}H}$$ and, by the
multiplicative property of $\Gamma$, we obtain
$$\mu_{A\times_{\beta}H}\co (\omega_{\alpha,\beta}\ot \omega_{\alpha,\beta})=p_{A\ot
H}\co\mu_{A\ot_{\beta}H}\co (\Gamma\ot \Gamma)\co (i_{A\ot H}\ot
i_{A\ot H})=p_{A\ot H}\co \Gamma \co \mu_{A\ot_{\alpha}H}\co
(i_{A\ot H}\ot i_{A\ot H})$$
 $$=\omega_{\alpha,\beta}\co
\mu_{A\times_{\alpha}H}.$$

Therefore, $\omega_{\alpha,\beta}$ is an isomorphism of algebras.

On the other hand, using (e1), (e2) and the property of left
$A$-module morphism of $\Gamma$ we have
$$\varphi_{A\times_{\beta}H}\co (A\ot \omega_{\alpha,\beta})=p_{A\ot
H}\co(\mu_{A}\ot H)\co (A\ot (\nabla_{A\ot H}\co\Gamma\co i_{A\ot
H}))=p_{A\ot H}\co(\mu_{A}\ot H)\co (A\ot (\Gamma\co i_{A\ot H}))$$
$$=p_{A\ot H}\co \Gamma\co (\mu_{A}\ot H)\co (A\ot  i_{A\ot H})=
p_{A\ot H}\co \Gamma\co \nabla_{A\ot H}\co(\mu_{A}\ot H)\co (A\ot
i_{A\ot H})=\omega_{\alpha,\beta}\co \varphi_{A\times_{\alpha}H}$$

and this proves that $\omega_{\alpha,\beta}$ is a morphism of left
$A$-modules.

Finally, using similar arguments and the property of right
$H$-comodule morphism of $\Gamma$ we obtain that
$\omega_{\alpha,\beta}$ is a morphism of right $H$-comodules
because:
$$\rho_{A\times_{\beta}H}\co \omega_{\alpha,\beta}=(p_{A\ot H}\ot H)\co
(A\ot \delta_{H})\co \nabla_{A\ot H}\co \Gamma \co i_{A\ot
H}=(p_{A\ot H}\ot H)\co\rho_{A\ot_{\beta} H}\co \Gamma \co i_{A\ot
H}$$
$$=(((p_{A\ot H}\co\Gamma)\ot H)\co  \rho_{A\ot_{\alpha} H}\co i_{A\ot
H}=(((p_{A\ot H}\co\Gamma\co \nabla_{A\ot H})\ot H)\co
\rho_{A\ot_{\alpha} H}\co i_{A\ot H}=(\omega_{\alpha,\beta}\ot H)\co
\rho_{A\times_{\alpha} H}.$$

\begin{rem}
{\rm By the previous theorem, we obtain that the notion of
equivalent crossed products is the one used in \cite{ana1} in a
category of modules over a commutative ring. Following the
terminology used in \cite{ana1}, the pair of morphisms $f_{\Gamma}$
and $f_{\Gamma}^{-1}$ is an example of gauge transformation. Also,
this notion is a generalization of the one that we can find in the
Hopf algebra world (see \cite{doi}, \cite{Guccione}).

The following results, Theorem \ref{principal-2} and Corollary
\ref{cor1-principal-2} will be used in Theorem \ref{teo39} to obtain
the meaning of the notion of equivalence between two weak crossed
products in terms of morphisms of $Reg_{\varphi_{A}}(H,A)$. Note
that this characterization it is the key to prove the main result of
this section, that is, Theorem \ref{principal}. }
\end{rem}

\begin{teo}
\label{principal-2} Let $\Gamma $ and $f_{\Gamma}$ as in Remark
\ref{relevant-remark} and such that
\begin{equation}
\label{prin-condition} \Gamma\co \nabla_{A\ot H}=\nabla_{A\ot
H}\co\Gamma=\Gamma.
\end{equation}
Under the hypothesis of Theorem \ref{principal-1}, $\Gamma $  is a
multiplicative morphism that preserves the preunit $\nu=\nabla_{A\ot
H}\co (\eta_{A}\ot \eta_{H})$ if and only if the following
equalities hold:
\begin{equation}
\label{preunit-crossed} p_{A\ot H}\co \Gamma \co \nu=p_{A\ot H}\co
\nu
\end{equation}
\begin{equation}
\label{psi-crossed} \mu_{A}\co (A\ot f_{\Gamma})\co
\psi_{H}^{A}=\mu_{A}\co (f_{\Gamma}\ot \varphi_{A})\co
(\delta_{H}\ot A)
\end{equation}
\begin{equation}
\label{sigma-crossed} \mu_{A}\co (A\ot f_{\Gamma})\co
\sigma_{H,\alpha}^{A}=\mu_{A}\co (\mu_{A}\ot \beta)\co (A\ot
\psi_{H}^{A}\ot H)\co (((f_{\Gamma}\ot H)\co \delta_{H})\ot
((f_{\Gamma}\ot H)\co \delta_{H}))
\end{equation}

Moreover, if $\Gamma$ preserves the preunit,  we have that
\begin{equation}
\label{f-gamma-eta} f_{\Gamma}\co \eta_{H}=\eta_{A}.
\end{equation}

\end{teo}

{\em Proof}: Assume that $\Gamma $  is a multiplicative morphism
that preserves the preunit. Then (\ref{preunit-crossed}) follows
easily and, by (\ref{prin-condition}), we have
\begin{equation}
\label{gamma-eta} \Gamma\co (A\ot \eta_{H})=\nabla_{A\ot H}\co (A\ot
\eta_{H})
\end{equation}
because
$$\Gamma\co (A\ot \eta_{H})=(\mu_{A}\ot H)\co (A\ot (\Gamma\co (\eta_{A}\ot
\eta_{H})))=(\mu_{A}\ot H) \co (A\ot (\Gamma\co \nu))=(\mu_{A}\ot
H)\co (A\ot
 \nu)$$
$$=\nabla_{A\ot H}\co (A\ot \eta_{H}).$$

On the other hand, the multiplicative condition for $\Gamma$ implies
that:
$$\Gamma \co \mu_{A\otimes_{\alpha}H}\co (\eta_{A}\ot H\ot A\ot
\eta_{H})=\mu_{A\otimes_{\beta}H}\co (\Gamma\ot \Gamma)\co
(\eta_{A}\ot H\ot A\ot \eta_{H}).$$ Equivalently
$$\Gamma \co (\mu_{A}\ot H)\co (A\ot (\sigma_{H,\alpha}^{A}\co (H\ot
\eta_{H})))\co \psi_{H}^{A}$$
\begin{equation}
\label{equ-1} =(\mu_{A}\ot H)\co (\mu_{A}\ot
\sigma_{H,\beta}^{A})\co (A\ot \psi_{H}^{A}\ot H)\co ((\Gamma\co
(\eta_{A}\ot H))\ot (\Gamma\co (A\ot \eta_{H})).
\end{equation}

By the normal condition for $\alpha$  we have
\begin{equation}
\label{sigma-alpha-1} \sigma_{H,\alpha}^{A}\co (H\ot
\eta_{H})=((\alpha\co (H\ot \Pi_{H}^{R})\co \delta_{H})\ot H)\co
\delta_{H}=(u_{1}\ot H)\co \delta_{H}=\nabla_{A\ot H}\co
(\eta_{A}\ot H) \end{equation} and then the upper side of
(\ref{equ-1}) is equal to $\Gamma\co \psi_{H}^{A}$. For the lower
side of (\ref{equ-1}) the following holds:

\begin{itemize}

\item[ ]$\hspace{0.38cm} (\mu_{A}\ot H)\co (\mu_{A}\ot
\sigma_{H,\beta}^{A})\co (A\ot \psi_{H}^{A}\ot H)\co ((\Gamma\co
(\eta_{A}\ot H))\ot (\Gamma\co (A\ot \eta_{H}))$

\item [ ]$= (\mu_{A}\ot H)\co (f_{\Gamma}\ot ((\mu_{A}\ot H)\co (A\ot \sigma_{H,\beta}^{A})\co
(\psi_{H}^{A}\ot H)\co (H\ot (\nabla_{A\ot H}\co (A\ot
\eta_{H})))))\co (\delta_{H}\ot A) $

\item [ ]$=(\mu_{A}\ot H)\co (f_{\Gamma}\ot ((\mu_{A}\ot H)\co (A\ot \sigma_{H,\beta}^{A})\co
(\psi_{H}^{A}\ot H)\co (H\ot A\ot \eta_{H})))\co (\delta_{H}\ot A) $

\item [ ]$=(\mu_{A}\ot H)\co (f_{\Gamma}\ot ((\mu_{A}\ot H)\co (A\ot (\nabla_{A\ot H}\co
(\eta_{A}\ot H)))))\co (H\ot \psi_{H}^{A})\co (\delta_{H}\ot A) $

\item [ ]$=(\mu_{A}\ot H)\co (f_{\Gamma}\ot \psi_{H}^{A})\co (\delta_{H}\ot A) $

\end{itemize}

where the first equality follows by (\ref{aux-import}) and
(\ref{gamma-eta}), the second one by (\ref{c11}), the third one by
(\ref{sigma-alpha-1}) and the fourth one by the properties of
$\nabla_{A\ot H}$.

Thus, (\ref{equ-1}) is equivalent to
\begin{equation}
\label{new-equ-1} \Gamma\co \psi_{H}^{A}=(\mu_{A}\ot H)\co
(f_{\Gamma}\ot \psi_{H}^{A})\co (\delta_{H}\ot A)
\end{equation}
and then composing in both sides with $A\ot \varepsilon_{H}$ we get
(\ref{psi-crossed}).

Also, the  multiplicative condition for $\Gamma$ implies the
following:
$$\Gamma \co \mu_{A\otimes_{\alpha}H}\co (\eta_{A}\ot H\ot \eta_{A}\ot
H)=\mu_{A\otimes_{\beta}H}\co (\Gamma\ot \Gamma)\co (\eta_{A}\ot
H\ot \eta_{A}\ot H).$$ Equivalently
$$\Gamma \co (\mu_{A}\ot H)\co (A\ot \sigma_{H,\alpha}^{A})\co
((\nabla_{A\ot H}\co (\eta_{A}\ot H))\ot H)$$
\begin{equation}
\label{equ-2} =(\mu_{A}\ot H)\co (\mu_{A}\ot
\sigma_{H,\beta}^{A})\co (A\ot \psi_{H}^{A}\ot H)\co ((\Gamma\co
(\eta_{A}\ot H))\ot (\Gamma\co (\eta_{A}\ot H)).
\end{equation}

Therefore, by (\ref{aw1}) and (\ref{aux-import}) we obtain that
(\ref{equ-2}) is equivalent to
\begin{equation}
\label{new-equ-2} \Gamma\co \sigma_{H,\alpha}^{A}=(\mu_{A}\ot H)\co
(\mu_{A}\ot \sigma_{H,\beta}^{A})\co (A\ot \psi_{H}^{A}\ot H)\co
((f_{\Gamma}\ot H)\co \delta_{H})\ot ((f_{\Gamma}\ot H)\co
\delta_{H})).
\end{equation}
Composing in both sides with $A\ot \varepsilon_{H}$ and using (iii)
of Proposition \ref{sigma-prop} we obtain (\ref{sigma-crossed}).

Conversely, assume that (\ref{preunit-crossed}), (\ref{psi-crossed})
and (\ref{sigma-crossed}) hold. Then,
$$\Gamma\co \nu=\nabla_{A\ot H}\co \Gamma \co \nu=\nabla_{A\ot H}\co
\nu=\nu$$ and $\Gamma$ preserves the preunit. Moreover, to prove
that $\Gamma$ is multiplicative first we show that, if
(\ref{psi-crossed}) holds, then (\ref{new-equ-1}) holds and
similarly for (\ref{sigma-crossed}) and (\ref{new-equ-2}). Indeed:

\begin{itemize}

\item[ ]$\hspace{0.38cm}\Gamma\co \psi_{H}^{A} $

\item [ ]$= ((\mu_{A}\co (A\ot f_{\Gamma}))\ot H)\co (A\ot
\delta_{H})\co \psi_{H}^{A}$

\item [ ]$= ((\mu_{A}\co (A\ot f_{\Gamma})\co \psi_{H}^{A})\ot H)\co
(H\ot c_{H,A})\co (\delta_{H}\ot A)  $

\item [ ]$=((\mu_{A}\co (f_{\Gamma}\ot \varphi_{A})\co
(\delta_{H}\ot A))\ot H)\co (H\ot c_{H,A})\co (\delta_{H}\ot A) $

\item [ ]$=(\mu_{A}\ot H)\co
(f_{\Gamma}\ot \psi_{H}^{A})\co (\delta_{H}\ot A) $

\end{itemize}

The first equality follows by (\ref{f-gamma}), the second and the
last ones by the coassociativity of $\delta_{H}$ and the third one
by (\ref{psi-crossed}).

\begin{itemize}

\item[ ]$\hspace{0.38cm}\Gamma\co \sigma_{H,\alpha}^{A} $

\item [ ]$=((\mu_{A}\co (A\ot f_{\Gamma}))\ot H)\co (A\ot
\delta_{H})\co \sigma_{H,\alpha}^{A} $

\item [ ]$= ((\mu_{A}\co (A\ot f_{\Gamma})\co
\sigma_{H,\alpha}^{A})\ot \mu_{H})\co \delta_{H^{2}}$

\item [ ]$= ((\mu_{A}\co (\mu_{A}\ot \beta)\co (A\ot
\psi_{H}^{A}\ot H)\co (((f_{\Gamma}\ot H)\co \delta_{H})\ot
((f_{\Gamma}\ot H)\co \delta_{H})))\ot \mu_{H})\co \delta_{H^{2}}$

\item [ ]$=(\mu_{A}\ot H)\co
(\mu_{A}\ot \sigma_{H,\beta}^{A})\co (A\ot \psi_{H}^{A}\ot H)\co
(((f_{\Gamma}\ot H)\co \delta_{H})\ot ((f_{\Gamma}\ot H)\co
\delta_{H})) $

\end{itemize}

The first equality follows by (\ref{f-gamma}), the second one by
(\ref{delta-sigmaHA}), the third one by (\ref{sigma-crossed}) and
the last one by the definition of $\psi_{H}^{A}$, the naturality of
$c$ and the coassociativity of $\delta_{H}$.

Then,

\begin{itemize}

\item[ ]$\hspace{0.38cm} \Gamma \co \mu_{A\ot_{\alpha}H}$

\item [ ]$=((\mu_{A}\co (A\ot f_{\Gamma}))\ot H)\co (\mu_{A}\ot \delta_{H})\co
(\mu_{A}\ot \sigma_{H,\alpha}^{A})\co (A\ot \psi_{H}^{A}\ot H) $

\item [ ]$=(\mu_{A}\ot H)\co (\mu_{A}\ot (\Gamma\co \sigma_{H,\alpha}^{A}))\co
(A\ot \psi_{H}^{A}\ot H) $

\item [ ]$= (\mu_{A}\ot H)\co (\mu_{A}\ot \sigma_{H,\beta}^{A})\co (\mu_{A}\ot ((\mu_{A}\ot H)\co
(f_{\Gamma}\ot \psi_{H}^{A})\co (\delta_{H}\ot A))\ot H)\co$
\item[ ]$\hspace{0.38cm}(A\ot \psi_{H}^{A}\ot ((f_{\Gamma}\ot H)\co \delta_{H}))  $

\item [ ]$=(\mu_{A}\ot H)\co (\mu_{A}\ot \sigma_{H,\beta}^{A})\co (A\ot (\Gamma \co
((\mu_A\ot H)\co (A\ot \psi_{H}^{A})\co (\psi_{H}^{A}\ot A)))\ot
H)\co$
\item[ ]$\hspace{0.38cm}(A\ot H\ot A\ot ((f_{\Gamma}\ot H)\co \delta_{H}))  $

\item [ ]$= (\mu_{A}\ot H)\co (A\ot \mu_{A}\ot H)\co (A\ot A\ot \sigma_{H,\beta}^{A})\co
(A\ot ( \Gamma\co \psi_{H}^{A})\ot H)\co (A\ot H\ot \Gamma) $

\item [ ]$= (\mu_{A}\ot H)\co (A\ot \mu_{A}\ot H)\co (A\ot A\ot \sigma_{H,\beta}^{A})\co
(A\ot ( (\mu_{A}\ot H)\co (f_{\Gamma}\ot \psi_{H}^{A})\co
(\delta_{H}\ot A))\ot H)\co $
\item[ ]$\hspace{0.38cm}(A\ot H\ot   \Gamma) $

\item [ ]$=\mu_{A\ot_{\beta}H} \co (\Gamma\ot \Gamma) $

\end{itemize}

The first equality follows by (\ref{f-gamma}), the second and the
last ones by the associativity of $\mu_{A}$, the third one by
(\ref{new-equ-2}), the fourth and the sixth ones by
(\ref{new-equ-1}) and the left $A$-linearity of $\Gamma$, the fifth
one by (\ref{wmeas-wcp}).

Finally, (\ref{f-gamma-eta}) follows by:
$$f_{\Gamma}\co \eta_{H}=(A\ot \varepsilon_{H})\co \Gamma\co
(\eta_{A}\ot \eta_{H})=(A\ot \varepsilon_{H})\co \Gamma\co
\nabla_{A\ot H}\co (\eta_{A}\ot \eta_{H})$$
$$=(A\ot \varepsilon_{H})\co
\Gamma\co \nu=(A\ot \varepsilon_{H})\co \nu=\eta_{A}.$$

\begin{cor}
\label{cor1-principal-2} Under the hypothesis of Theorem
\ref{principal-2}, if (\ref{psi-crossed}) holds,
(\ref{sigma-crossed}) is equivalent to
\begin{equation}
\label{sigma-crossed-2} \mu_{A}\co (A\ot f_{\Gamma})\co
\sigma_{H,\alpha}^{A}=[\mu_{A}\co ((\varphi_{A}\co (H\ot
f_{\Gamma}))\ot f_{\Gamma})\co (H\ot c_{H,H})\co (\delta_{H}\ot
H)]\wedge \beta.
\end{equation}

Then, if $f_{\Gamma}\in Reg_{\varphi_{A}}(H,A)$, we obtain that
(\ref{sigma-crossed}) is equivalent to
\begin{equation}
\label{cohomologous-1} \alpha\wedge \partial_{1,1}(f_{\Gamma})=
\partial_{1,0}(f_{\Gamma})\wedge \partial_{1,2}(f_{\Gamma}) \wedge
\beta.
\end{equation}

\end{cor}

{\em Proof}: If (\ref{sigma-crossed-2}) holds, by
(\ref{psi-crossed}), the naturality of $c$ and the coassociativity
of $\delta_{H}$, we obtain (\ref{sigma-crossed}):

\begin{itemize}

\item[ ]$\hspace{0.38cm}\mu_{A}\co (A\ot f_{\Gamma})\co
\sigma_{H,\alpha}^{A}$

\item [ ]$=[\mu_{A}\co ((\varphi_{A}\co (H\ot
f_{\Gamma}))\ot f_{\Gamma})\co (H\ot c_{H,H})\co (\delta_{H}\ot
H)]\wedge \beta $

\item [ ]$= \mu_{A}\co ((\mu_{A}\co (A\ot f_{\Gamma})\co \psi_{H}^{A})\ot \beta)\co
 (H\ot c_{H,A}\ot H)\co (\delta_{H}\ot ((f_{\Gamma}\ot H)\co \delta_{H}))$

\item [ ]$=\mu_{A}\co ((\mu_{A}\co (f_{\Gamma}\ot \varphi_{A})\co (\delta_{H}\ot A))\ot \beta)\co
 (H\ot c_{H,A}\ot H)\co (\delta_{H}\ot ((f_{\Gamma}\ot H)\co \delta_{H})) $

\item [ ]$=\mu_{A}\co (\mu_{A}\ot \beta)\co (A\ot
\psi_{H}^{A}\ot H)\co (((f_{\Gamma}\ot H)\co \delta_{H})\ot
((f_{\Gamma}\ot H)\co \delta_{H})) $

\end{itemize}

On the other hand, (\ref{psi-crossed}) holds we have that
(\ref{new-equ-1}) holds and then if we assume (\ref{sigma-crossed}),
using (\ref{f-gamma}), the definition of $\psi_{H}^{A}$, the
naturality of $c$ and the coassociativity of $\delta_{H}$, we
obtain:

\begin{itemize}

\item[ ]$\hspace{0.38cm}\mu_{A}\co (A\ot f_{\Gamma})\co
\sigma_{H,\alpha}^{A}$

\item [ ]$= \mu_{A}\co (\mu_{A}\ot \beta)\co (A\ot \psi_{H}^{A}\ot H)\co
(((f_{\Gamma}\ot H)\co \delta_{H})\ot ((f_{\Gamma}\ot H)\co
\delta_{H}))$

\item [ ]$= \mu_{A}\co (A\ot \beta)\co ((\Gamma\co \psi_{H}^{A})\ot H)\co
(H\ot ((f_{\Gamma}\ot H)\co \delta_{H}))$

\item [ ]$=\mu_{A}\co (A\ot \beta)\co
(((((\mu_{A}\co (A\ot f_{\Gamma}))\ot H)\co (A\ot \delta_{H})) \co
\psi_{H}^{A})\ot H)\co (H\ot ((f_{\Gamma}\ot H)\co \delta_{H})) $

\item [ ]$=[\mu_{A}\co ((\varphi_{A}\co (H\ot
f_{\Gamma}))\ot f_{\Gamma})\co (H\ot c_{H,H})\co (\delta_{H}\ot
H)]\wedge \beta $

\end{itemize}

Finally, it is obvious that
\begin{equation}
\label{cohomologous-11} \mu_{A}\co (A\ot f_{\Gamma})\co
\sigma_{H,\alpha}^{A}=\alpha\wedge \partial_{1,1}(f_{\Gamma})
\end{equation}
and,  by (\ref{firts-cocom-1}) and $\beta\co \Omega_{H}^{2}=\beta$
we have
$$\partial_{1,0}(f_{\Gamma})\wedge
\partial_{1,2}(f_{\Gamma}) \wedge \beta=$$
\begin{equation}
\label{cohomologous-12} [\mu_{A}\co ((\varphi_{A}\co (H\ot
f_{\Gamma}))\ot f_{\Gamma})\co (H\ot c_{H,H})\co (\delta_{H}\ot
H)]\wedge \beta.
\end{equation}

\begin{teo}
\label{teo39} Under the hypothesis of Theorem \ref{principal-1}, the
weak crossed products $A\ot_{\alpha} H$, $A\ot_{\beta} H$
associated to $\alpha$ and $\beta$ are equivalent if and only if
there exists $f\in Reg_{\varphi_{A}}^{+}(H,A)$ such that the
equalities (\ref{psi-crossed}) and (\ref{cohomologous-1}) hold.
\end{teo}

{\em Proof}: If the weak crossed products $A\ot_{\alpha} H$,
$A\ot_{\beta} H$ are equivalent, by Theorem \ref{principal-1}, there
exist  multiplicative and preunit preserving morphisms $\Gamma,
\Gamma^{\prime} \in\;_{A}Hom^{H}_{\mathcal C}(A\ot H, A\ot H)$
satisfying (e1)-(e3). Then, by Remark \ref{relevant-remark},
$f_{\Gamma}\in Reg_{\varphi_{A}}(H,A)$, and by Theorem
\ref{principal-2}, the equalities (\ref{psi-crossed}) and
$f_{\Gamma}\co \eta_{H}=\eta_{A}$ hold. Finally, by Corollary
\ref{cor1-principal-2} we get (\ref{cohomologous-1}). Conversely,
let $f\in Reg_{\varphi_{A}}^{+}(H,A)$, with inverse $f^{-1}$. Then,
 $\Gamma_{f}$ and $\Gamma_{f^{-1}}$ are morphisms of left
$A$-modules and right $H$-comodules satisfying (e1)-(e3)  and
preserving the preunit $\nu=\nabla_{A\ot H}\co (\eta_{A}\ot
\eta_{H})$. Indeed: By (\ref{delta-pi31}) and (iii) of Proposition
\ref{f-eta}, we have
$$\Gamma_{f}\co \nu=(f\ot H)\co \delta_{H}\co \eta_{H}= ((f\co
\overline{\Pi}_{H}^{L})\ot H)\co \delta_{H}\co \eta_{H}=(u_{1}\ot
H)\co \delta_{H}\co \eta_{H}=\nu.$$ Similarly, $\Gamma_{f^{-1}}\co
\nu=\nu$. By Theorem \ref{principal-2} and Corollary
\ref{cor1-principal-2}, $\Gamma_{f}$ is multiplicative and
$$\omega_{\alpha,\beta}=p_{A\ot H}\co \Gamma_{f}\co i_{A\ot H}$$
is an $H$-comodule algebra isomorphism with inverse
$\omega_{\alpha,\beta}^{-1}=p_{A\ot H}\co \Gamma_{f^{-1}}\co i_{A\ot
H}$. Then, $\Gamma_{f^{-1}}$ is multiplicative and, by Theorem
\ref{principal-1}, we obtain that $A\ot_{\alpha} H$, $A\ot_{\beta}
H$ are equivalent.

\begin{rem}
\label{second-relevant-remark}
 {\rm Note that, by the previous Theorem, if $A\ot_{\alpha} H$, $A\ot_{\beta} H$
 are equivalent and $f\in
Reg_{\varphi_{A}}^{+}(H,A)$ is the morphism inducing the
equivalence, by Theorem \ref{principal-2}, and Corollary
\ref{cor1-principal-2} we also have
\begin{equation}
\label{psi-crossed-f-1} \mu_{A}\co (A\ot f^{-1})\co
\psi_{H}^{A}=\mu_{A}\co (f^{-1}\ot \varphi_{A})\co (\delta_{H}\ot
A),
\end{equation}
\begin{equation}
\label{new-equ-1-f-1} \Gamma_{f^{-1}}\co \psi_{H}^{A}=(\mu_{A}\ot
H)\co (f^{-1}\ot \psi_{H}^{A})\co (\delta_{H}\ot A)
\end{equation}
\begin{equation}
\label{sigma-crossed-f-1} \mu_{A}\co (A\ot f^{-1})\co
\sigma_{H,\beta}^{A}=\mu_{A}\co (\mu_{A}\ot \alpha)\co (A\ot
\psi_{H}^{A}\ot H)\co (((f^{-1}\ot H)\co \delta_{H})\ot ((f^{-1}\ot
H)\co \delta_{H})),
\end{equation}
\begin{equation}
\label{new-equ-2-f-1} \Gamma_{f^{-1}}\co
\sigma_{H,\beta}^{A}=(\mu_{A}\ot H)\co (\mu_{A}\ot
\sigma_{H,\alpha}^{A})\co (A\ot \psi_{H}^{A}\ot H)\co ((f^{-1}\ot
H)\co \delta_{H})\ot ((f^{-1}\ot H)\co \delta_{H})),
\end{equation}
\begin{equation}
\label{sigma-crossed-2-f-1} \mu_{A}\co (A\ot f^{-1})\co
\sigma_{H,\beta}^{A}=[\mu_{A}\co ((\varphi_{A}\co (H\ot f^{-1}))\ot
f^{-1})\co (H\ot c_{H,H})\co (\delta_{H}\ot H)]\wedge \alpha,
\end{equation}
and
\begin{equation}
\label{cohomologous-1-f-1} \beta\wedge \partial_{1,1}(f^{-1})=
\partial_{1,0}(f^{-1})\wedge \partial_{1,2}(f^{-1}) \wedge
\alpha.
\end{equation}

}
\end{rem}

\begin{rem}
\label{third-relevant-remark}
 {\rm Note that, if $H$ is a cocommutative weak Hopf algebra, $(A,\varphi_{A})$
is a weak left $H$-module algebra and $f:H\rightarrow A$ is a
morphism, the equality (\ref{psi-crossed}) it is always true if $A$
is commutative, because:
$$\mu_{A}\co (A\ot f)\co
\psi_{H}^{A}=\mu_{A}\co (\varphi_{A}\ot f)\co (H\ot c_{H,A})\co
((c_{H,H}\co \delta_{H})\ot A)= \mu_{A}\co c_{A,A}\co (f\ot
\varphi_{A})\co (\delta_{H}\ot A)$$
$$=\mu_{A}\co (f\ot
\varphi_{A})\co (\delta_{H}\ot A).$$

Then, if $(A,\varphi_{A})$ is a commutative left $H$-module algebra,
the equivalence between two weak crossed products $A\ot_{\alpha} H$,
$A\ot_{\beta} H$ is determined by the inclusion of $f$ in
$Reg_{\varphi_{A}}^{+}(H,A)$ and the equality
(\ref{cohomologous-1}). In this case (\ref{cohomologous-1}) is
equivalent to say that $\alpha\wedge \beta^{-1}\in
Im(D^{1+}_{\varphi_{A}})$.

 }
 \end{rem}

\begin{teo}
\label{principal} Let $H$ be a cocommutative weak Hopf algebra and
$(A,\varphi_{A})$ a commutative left $H$-module algebra. Then there
is a bijective correspondence between $H^{2}_{\varphi_{A}}(H,A)$ and
the equivalence classes of weak crossed products of $A\ot_{\alpha}
H$ where $\alpha: H\ot H\rightarrow A$ satisfy the 2-cocycle
condition (\ref{2-cocycle-sigma})(equivalently
(\ref{cocycle-equivalent})) and the normal condition
(\ref{normal-sigma}).
\end{teo}

{\em Proof}: First note that $H^{2}_{\varphi_{A}}(H,A)$ is
isomorphic to $H^{2+}_{\varphi_{A}}(H,A)$. Then, it is suffices to
prove the result for $H^{2+}_{\varphi_{A}}(H,A)$. Let $\alpha,
\beta\in Reg_{\varphi_{A}}^{+}(H^2,A)$ such that satisfies the
2-cocycle condition (\ref{2-cocycle-sigma}) (in the commutative case
the twisted condition it is always true). If $A\ot_{\alpha} H$,
$A\ot_{\beta} H$ are equivalent, by the previous remark, we have
that there exists $f$ in $Reg_{\varphi_{A}}^{+}(H,A)$ such  that
$\alpha\wedge \beta^{-1}\in Im(D^{1+}_{\varphi_{A}})$. Then,
$\alpha$ and $\beta$ are in the same class in
$H^{2+}_{\varphi_{A}}(H,A)$. Conversely, if $[\alpha]=[\beta]$ in
$H^{2+}_{\varphi_{A}}(H,A)$, $\alpha$ and $\beta$ satisfies
(\ref{cohomologous-1}), i.e. $\alpha\wedge
\beta^{-1}=D^{1+}_{\varphi_{A}}(f)$, for $f \in
Reg_{\varphi_{A}}^{+}(H,A)$. Then, if $\Gamma_{f}$ is the morphism
defined in Remark \ref{relevant-remark}, we have that $\Gamma_{f}$
satisfies (\ref{preunit-crossed}), because, using that $f \in
Reg_{\varphi_{A}}^{+}(H,A)$, we obtain
$$p_{A\ot H}\co \Gamma_{f}\co \nu=p_{A\ot H}\co (f\ot H)\co
\delta_{H}\co \eta_{H}=p_{A\ot H}\co ((f\co
\overline{\Pi}_{H}^{L})\ot H)\co \delta_{H}\co \eta_{H}$$
$$=p_{A\ot
H}\co ((f\co \Pi_{H}^{L})\ot H)\co \delta_{H}\co \eta_{H}=p_{A\ot
H}\co (u_{1}\ot H)\co \delta_{H}\co \eta_{H}=p_{A\ot H}\co \nu.$$

In a similar way, $\beta\wedge
\alpha^{-1}=D^{1+}_{\varphi_{A}}(f^{-1})$ and $\Gamma_{f^{-1}}$
satisfies (\ref{preunit-crossed}). Then, by Theorem
\ref{principal-2}, $\Gamma_{f}$ and $\Gamma_{f^{-1}}$ are
multiplicative morphisms of left $A$-modules and right $H$-comodules
preserving the preunit and satisfying (e1)-(e3). Therefore, by
Theorem \ref{principal-1}, we obtain that $A\ot_{\alpha} H$,
$A\ot_{\beta} H$  are equivalent weak crossed products.

\section*{Acknowledgements}
The authors were supported by  Ministerio de Ciencia e Innovaci\'on
(Project: MTM2010-15634) and by FEDER.

\end{document}